\documentclass{amsart}
\usepackage{amssymb}
\newtheorem{theorem}{Theorem}[section]
\newtheorem{lemma}[theorem]{Lemma}
\newtheorem{proposition}[theorem]{Proposition}

\newtheorem{definition}[theorem]{Definition}
\newtheorem{conj}[theorem]{Conjecture}
\newtheorem{remark}[theorem]{Remark}

\newcommand{\C}{\mathbb{C}}
\newcommand{\Z}{\mathbb{Z}}
\newcommand{\CP}{\mathbb{CP}}
\numberwithin{equation}{section}
\title[Fundamental groups and symplectic invariants]{Fundamental groups of 
complements of plane curves and symplectic invariants}
\author{D. Auroux}
\address{Department of Mathematics, M.I.T., Cambridge MA 02139, USA}
\email{auroux@math.mit.edu}
\author{S.\,K. Donaldson}
\address{Department of Mathematics, Imperial College, London SW7 2BZ, 
United Kingdom}
\email{s.donaldson@ic.ac.uk}
\author{L. Katzarkov}
\address{Department of Mathematics, University of California, 
Irvine, CA 92697, USA}
\email{lkatzark@math.uci.edu}
\author{M. Yotov}
\address{Department of Mathematics, University of California, 
Irvine, CA 92697, USA}
\email{myotov@math.uci.edu}

\begin{document}
\begin{abstract}
Introducing the notion of stabilized fundamental group for the complement
of a branch curve in $\CP^2$, we define effectively computable 
invariants of symplectic 4-manifolds that generalize those previously
introduced by Moishezon and Teicher for complex projective surfaces.
Moreover, we study the structure of these invariants and formulate
conjectures supported by calculations on new examples.
\end{abstract}

\maketitle

\section{Introduction}

Using approximately holomorphic techniques first introduced in
\cite{D1}, it was shown in \cite{AK} (see also \cite{A})
that compact symplectic $4$-manifolds with integral
symplectic class can be realized as branched covers of $\CP^2$
and can be investigated using the braid group techniques developed by 
Moishezon and subsequently by
Moishezon and Teicher for the study of complex surfaces (see e.g.\
\cite{Tsurv}):

\begin{theorem}[\cite{AK}]
Let $(X,\omega)$ be a compact symplectic $4$-manifold, and let $L$ be
a line bundle with $c_1(L)=\frac{1}{2\pi}[\omega]$. Then
there exist branched covering maps $f_k:X\to\CP^2$ defined by
approximately holomorphic sections of $L^{\otimes k}$ for all large enough
values of $k$; the corresponding branch curves $D_k\subset\CP^2$ 
admit only nodes (both orientations) and complex cusps as singularities, 
and give rise to well-defined braid monodromy invariants. Moreover, up to
admissible creations and cancellations of pairs of nodes in the branch 
curve, for large $k$ the topology of $f_k$ is a symplectic invariant.
\end{theorem}

This makes it possible to associate to $(X,\omega)$ a sequence of invariants
(indexed by $k\gg 0$) consisting of two objects: the braid 
monodromy characterizing the branch curve $D_k$, and the {\it geometric 
monodromy representation} $\theta_k:\pi_1(\CP^2-D_k)\to S_n$ 
($n=\deg f_k$) characterizing the $n$-fold covering of $\CP^2-D_k$
induced by $f_k$ \cite{AK}. These invariants are extremely powerful (from
them one can recover $(X,\omega)$ up to symplectomorphism) but
too complicated to handle in practical cases. 

In the study of complex surfaces, Moishezon and Teicher have 
shown that the fundamental
group $\pi_1(\CP^2-D)$ (or, restricting to an affine subset,
$\pi_1(\C^2-D)$) can be computed explicitly in some simple examples; 
generally speaking, this group has 
been expected to provide a valuable invariant for distinguishing
diffeomorphism types of complex surfaces of general type. 
However, in the symplectic case, it is affected 
by creations and cancellations of pairs of nodes and cannot be used 
immediately as an invariant.

We will introduce in \S 2 a certain quotient $G_k$ (resp.\ $\bar{G}_k$) of 
$\pi_1(\C^2-D_k)$ (resp.\ $\pi_1(\CP^2-D_k)$), the
{\em stabilized fundamental group}, which remains invariant under creations
and cancellations of pairs of nodes. As an immediate corollary of the 
construction and of Theorem 1.1, we obtain the following

\begin{theorem}
For large enough $k$, the stabilized groups $G_k=G_k(X,\omega)$ 
$($resp.\ $\bar{G}_k(X,\omega))$ and their 
reduced subgroups $G_k^0=G_k^0(X,\omega)$
are symplectic invariants of the manifold $(X,\omega)$.
\end{theorem}

These invariants can be computed explicitly in various examples, some
due to Moishezon, Teicher and Robb, others new; these examples will be 
presented in \S 4, and a brief overview of the techniques involved in the 
computations is given in \S 6 and \S 7. The new examples include
double covers of $\CP^1\times\CP^1$ branched along arbitrary complex 
curves (Theorem 4.6 and \S 7)~; similar methods should apply to
other double covers as well, thus providing results for both types of 
so-called Horikawa surfaces.

The available data suggest several conjectures
about the structure of the stabilized fundamental groups. 

First of all, it appears that in most examples the stabilization operation
does not actually affect the fundamental group. The only known exceptions 
are given by ``small'' linear systems with insufficient ampleness 
properties, where the stabilization is a quotient by a non-trivial
subgroup (see \S 4). Therefore we have the following

\begin{conj}
Assume that $(X,\omega)$ is a complex surface, and let $D_k$ be
the branch curve of a generic projection to $\CP^2$ of the
projective embedding of $X$ given by the linear system $|kL|$. 
Then, provided that $k$ is large enough, the stabilization operation is
trivial, i.e.\ $G_k(X,\omega)\simeq \pi_1(\C^2-D_k)$ and
$\bar{G}_k(X,\omega)\simeq \pi_1(\CP^2-D_k)$.
\end{conj}

An important class of fundamental groups for which the conjecture holds
will be described in \S 3.

Moreover, the structure of the stabilized fundamental groups
seems to be remarkably simple, at least when the manifold $X$ is simply 
connected; in all known examples they are extensions of a symmetric group
by a solvable group, while there exist plane curves with much more
complicated complements \cite{DOZ,KP}.
In fact these groups seem to be largely determined by 
intersection pairing data in $H_2(X,\Z)$. More precisely,
the following result will be proved in \S 5:

\begin{definition}
Let $\Lambda_k$ be the image of the map $\lambda_k:H_2(X,\Z)\to\Z^2$ defined 
by $\lambda_k(\alpha)=(\alpha\cdot L_k,\alpha\cdot R_k)$, where 
$L_k=k\,c_1(L)$ and $R_k=c_1(K_X)+3L_k$ are the classes in $H^2(X,\Z)$ 
Poincar\'e dual to a hyperplane section and to the ramification curve 
respectively.
\end{definition}

\begin{theorem}
If the symplectic manifold $X$ is simply connected,
then there exists a natural surjective homomorphism
$\phi_k:\mathrm{Ab}\,G_k^0(X,\omega)\to (\Z^2/\Lambda_k)\otimes
\mathcal{R}_{n_k}\simeq (\Z^2/\Lambda_k)^{n_k-1}$, where 
$n_k=\deg f_k=L_k\cdot L_k$, and $\mathcal{R}_{n_k}$ is the reduced
regular representation of $S_{n_k}$ $($isomorphic to $\Z^{n_k-1})$.
\end{theorem}

The map $\phi_k$ is $(G_k,S_{n_k})$-equivariant, 
in the sense that $\phi_k(g^{-1}\gamma g)=\theta_k(g)\cdot \phi_k(\gamma)$ 
for any elements $g\in G_k(X,\omega)$ and 
$\gamma\in\mathrm{Ab}\,G_k^0(X,\omega)$
(cf.\ also Lemma 5.2).

In the examples discussed in \S 4, the group $G_k^0$ is always
close to being abelian, and $\phi_k$ is always an isomorphism. It seems
likely that the injectivity of $\phi_k$ can be proved using 
techniques similar to those described in \S 6--7.
Therefore, it makes sense to formulate the following

\begin{conj}
If the symplectic manifold $X$ is simply connected and $k$ is large
enough, then $\mathrm{Ab}\,G_k^0(X,\omega)\simeq (\Z^2/\Lambda_k)\otimes
\mathcal{R}_{n_k}$,
and the commutator subgroup $[G_k^0,G_k^0]$ is a quotient of $(\Z_2)^2$.
\end{conj}

Conjectures 1.3 and 1.6 provide an almost complete tentative description of the
structure of fundamental groups of branch curve complements in high degrees.
In relation with the property $(*)$ introduced in \S 3, they also provide a 
framework to explain various observations and conjectures made in 
\cite{Tstruct} and \cite{Robb2}.

Conjectures 1.3 and 1.6 seem to indicate that fundamental groups of branch
curve complements fail to provide useful invariants to symplectically 
distinguish homeomorphic manifolds. This is in sharp contrast with the 
braid monodromy data, which completely determines the symplectomorphism type
of $(X,\omega)$ \cite{AK}~; how to introduce effectively computable 
invariants retaining more of the information contained in the braid 
monodromy remains an open question.

\section{Braid monodromy and stabilized fundamental groups}

Let $D_k$ be the branch curve of a covering map $f_k:X\to\CP^2$
as in Theorem 1.1. Braid monodromy invariants are defined by considering
a generic projection $\pi:\CP^2-\{\mathrm{pt}\}\to\CP^1$: the pole of
the projection lies away from $D_k$, and a generic fiber of $\pi$ intersects
$D_k$ in $d=\deg D_k$ distinct points, the only exceptions being fibers
through cusps or nodes of $D_k$, or fibers that are tangent to $D_k$ at
one of its smooth points (``vertical tangencies''). Moreover we can assume
that the special points (cusps, nodes and vertical tangencies) of $D_k$ 
all lie in different fibers of $\pi$.

By restricting ourselves to an affine subset $\C^2\subset\CP^2$, choosing a
base point and trivializing the fibration $\pi$, we can view the monodromy 
of $\pi_{|D_k}$ as a group homomorphism from $\pi_1(\C-\{q_i\})$ (where $q_i$
are the images by $\pi$ of the special points of $D_k$) to the braid group 
$B_d$. More precisely, the monodromy around a vertical tangency is a 
half-twist (a braid that exchanges two of the $d$ intersection points of
the fiber with $D_k$ by rotating them around each other counterclockwise 
along a certain path); the monodromy 
around a positive (resp.\ negative) node is the square (resp.\ the inverse
of the square) of a half-twist; the monodromy around a cusp is
the cube of a half-twist \cite{Tsurv,AK}.

It is sometimes convenient to choose an ordered system of generating
loops for $\pi_1(\C-\{q_i\})$ (one loop going around each $q_i$), and
to express the monodromy as a {\em braid factorization}, i.e.\ a
decomposition of the central braid $\Delta^2$ (the monodromy around the
fiber at infinity, due to the non-triviality of the fibration $\pi$ over
$\CP^1$) into the product of the monodromies along the chosen generating
loops. However, this braid factorization is only well-defined up to 
simultaneous conjugation of all factors (i.e., a change in the
choice of the identification of the fibers with $\mathbb{R}^2$) and {\em
Hurwitz equivalence} (i.e., a rearrangement of the factors due to
a different choice of the system of generating loops).

The braid monodromy determines in a very explicit manner the fundamental
groups $\pi_1(\C^2-D_k)$ and $\pi_1(\CP^2-D_k)$. Indeed, consider a
generic fiber $\ell \simeq\C\subset\CP^2$ of the projection $\pi$ 
(e.g.\ the fiber containing the base point), intersecting
$D_k$ in $d$ distinct points. The free group $\pi_1(\ell-(\ell\cap
D_k))=F_d$ is generated by a system of $d$ loops
going around the various points in $\ell\cap D_k$.
The inclusion map 
$i:\ell-(\ell\cap D_k)\to \C^2-D_k$ induces a surjective 
homomorphism $i_*:F_d\to\pi_1(\C^2-D_k)$. 

\begin{definition}
The images of the standard generators of the free group $F_d$ and their 
conjugates are called {\it geometric generators} of 
$\pi_1(\C^2-D_k)$; the set of all geometric generators will be denoted by 
$\Gamma_k$.
\end{definition}

By the Zariski-Van Kampen theorem, $\pi_1(\C^2-D_k)$ is realized as a 
quotient of $F_d$ by relations corresponding to the various special points 
(vertical tangencies, nodes, cusps) of $D_k$; these relations express the
fact that the action of the braid monodromy on $F_d$ induces a trivial
action on $\pi_1(\C^2-D_k)$. To each factor in the
braid factorization one can associate a pair of elements $\gamma_1,\gamma_2\in
\Gamma_k$ (small loops around the two portions of $D_k$ that meet at the
special point), well-determined up to simultaneous conjugation.
The relation corresponding to a tangency is $\gamma_1\sim \gamma_2$;
for a node (of either orientation) it is $[\gamma_1,\gamma_2]\sim 1$; for a 
cusp it becomes $\gamma_1\gamma_2\gamma_1\sim\gamma_2\gamma_1\gamma_2$.
Taking into account all the special points of $D_k$ (i.e.\ considering the 
entire braid monodromy), we obtain a presentation of $\pi_1(\C^2-D_k)$.
Moreover, $\pi_1(\CP^2-D_k)$ is obtained from $\pi_1(\C^2-D_k)$ just by
adding the extra relation $g_1\dots g_d\sim 1$, where $g_i$ are the
images of the standard generators of $F_d$ under the inclusion.
\medskip

It follows from this discussion that the creation or cancellation of a
pair of nodes in $D_k$ may affect $\pi_1(\C^2-D_k)$ and $\pi_1(\CP^2-D_k)$
by adding or removing commutation relations between geometric generators.
Although it is reasonable to expect that negative nodes can always be
cancelled in the branch curves given by Theorem 1.1, the currently available
techniques are insufficient to prove such a statement. Instead, a more
promising approach is to compensate for these changes in the
fundamental groups by considering certain quotients where one stabilizes
the group by adding commutation relations between geometric generators.
The resulting group is in some sense more natural than $\pi_1(\C^2-D_k)$ 
from the symplectic point of view, and as a side benefit it is often easier 
to compute (see \S 7). Moreover, it also turns out that, in many cases, no 
information is lost in the stabilization process (see \S 3).

In order to define the stabilized group $G_k$, first observe that, because
the branching index of $f_k$ above a smooth point of $D_k$ is always $2$,
the geometric monodromy representation morphism 
$\theta_k:\pi_1(\CP^2-D_k)\to S_n$ describing the topology of the covering
above $\CP^2-D_k$ maps all geometric generators to 
transpositions in $S_n$. As seen above, to each nodal point of $D_k$ one 
can associate geometric generators $\gamma_1,\gamma_2\in\Gamma_k$, one
for each of the two intersecting portions of $D_k$, so that the
corresponding relation in $\pi_1(\C^2-D_k)$ is $[\gamma_1,\gamma_2]\sim 1$.
Since the branching
occurs in disjoint sheets of the cover, the two transpositions
$\theta_k(\gamma_1)$ and $\theta_k(\gamma_2)$ are necessarily disjoint
(i.e.\ they are distinct and commute). Therefore,
adding or removing pairs of nodes amounts to adding or removing relations
given by commutators of geometric generators associated to disjoint
transpositions.

\begin{definition}
Let $K_k$ $($resp.\ $\bar{K}_k)$ be the normal subgroup of $\pi_1(\C^2-D_k)$ 
$($resp.\ $\pi_1(\CP^2-D_k))$ generated by all commutators 
$[\gamma_1,\gamma_2]$ where $\gamma_1,\gamma_2\in\Gamma_k$ are such that 
$\theta_k(\gamma_1)$ and $\theta_k(\gamma_2)$ are disjoint transpositions.
The {\em stabilized fundamental group} is defined as
$G_k=\pi_1(\C^2-D_k)/K_k$, resp.\ $\bar{G}_k=\pi_1(\CP^2-D_k)/\bar{K}_k$.
\end{definition}

Certain natural subgroups of $G_k$ and $\bar{G}_k$ will play an important
role in the following sections. Define the {\it linking number}
homomorphism $\delta_k:\pi_1(\C^2-D_k)\to \Z$ by $\delta_k(\gamma)=1$ 
for every $\gamma\in\Gamma_k$; similarly one can define
$\bar\delta_k:\pi_1(\CP^2-D_k)\to \Z_d$. When $D_k$ is irreducible (which
is the general case), these can also be thought of as
abelianization maps from the fundamental groups to the homology groups
$H_1(\C^2-D_k,\Z)\simeq\Z$ and $H_1(\CP^2-D_k,\Z)\simeq \Z_d$.

\begin{lemma}
$\mathrm{Ker}\,\delta_k\simeq\mathrm{Ker}\,\bar\delta_k$.
\end{lemma}

\proof Since $\pi_1(\CP^2-D_k)=\pi_1(\C^2-D_k)/\langle g_1\dots g_d\rangle$
and $\delta_k(g_1\dots g_d)=d$, it is sufficient to show that the product
$g_1\dots g_d$
belongs to the center of $\pi_1(\C^2-D_k)$. Observe that the relation
in $\pi_1(\C^2-D_k)$ coming from a special point of $D_k$ can be rewritten
in the form $g\sim b_* g$ $\forall g\in F_d$, where $b\in
B_d$ is the braid monodromy around the given special point, acting on $F_d$.
In particular, if we consider the braid monodromy as a factorization 
$\Delta^2=\prod b_i$, we obtain that $g\sim (\prod b_i)_* g=(\Delta^2)_* g$
for any element $g$. However the action of the braid $\Delta^2$ on $F_d$
is exactly conjugation by $g_1\dots g_d$; we conclude that $g_1\dots g_d$
commutes with any element of $\pi_1(\C^2-D_k)$, hence the result.
\endproof

The homomorphisms $\delta_k$ and $\bar\delta_k$ are obviously surjective.
Moreover, $\theta_k$ is also surjective, because of the connectedness of 
$X$: the subgroup $\mathrm{Im}\,\theta_k\subseteq S_n$ is generated by 
transpositions and acts transitively on $\{1,\dots,n\}$, so it is equal 
to $S_n$. However, the image of 
$\theta^+_k=(\theta_k,\delta_k):\pi_1(\C^2-D_k)\to S_n\times\Z$
is the index 2 subgroup
$\{(\sigma,i): \mathrm{sgn}(\sigma)\equiv i\,\mathrm{mod}\,2\}$,
and similarly for $\bar\theta^+_k=(\theta_k,\bar\delta_k):
\pi_1(\CP^2-D_k)\to S_n\times\Z_d$ (note
that $d$ is always even).
Since $K_k\subseteq \mathrm{Ker}\,\theta^+_k$,
we can make the following definition:

\begin{definition}
Let $H_k^0=\mathrm{Ker}\,\theta^+_k\simeq \mathrm{Ker}\,\bar\theta^+_k$.
The {\em reduced subgroup} of $G_k$ is $G_k^0=H_k^0/K_k$.
We have the following exact sequences:
$$1\longrightarrow G_k^0\longrightarrow G_k\longrightarrow S_n\times\Z
\longrightarrow \Z_2 \longrightarrow 1,$$
$$1\longrightarrow G_k^0\longrightarrow \bar{G}_k\longrightarrow 
S_n\times\Z_d \longrightarrow \Z_2 \longrightarrow 1.$$
\end{definition}

Theorem 1.2 is now obvious from the definitions and from Theorem 1.1:
since creating a pair of nodes amounts to adding a relation of the form 
$[\gamma_1,\gamma_2]\sim 1$ where $[\gamma_1,\gamma_2]\in K_k$ (resp.\
$\bar{K}_k$), by construction it does not affect the groups
$G_k$, $\bar{G}_k$ and $G_k^0$, which are 
therefore symplectic invariants for $k$ large enough.

\section{$\tilde{B}_n$-groups and their stabilizations}

Denote by $B_n$ (resp.\ $P_n$, $P_{n,0}$) the braid group on $n$ strings
(resp.\ the subgroups of pure braids and pure braids of degree $0$), and
denote by $X_1,\dots,X_{n-1}$ the standard generators of $B_n$. Recall that
$X_i$ is a half-twist along a segment joining the points $i$ and $i+1$, and
that the relations among these generators are $[X_i,X_j]=1$ if $|i-j|\ge 2$
and $X_iX_{i+1}X_i=X_{i+1}X_iX_{i+1}$.

Let $\tilde{B}_n$ be the quotient of $B_n$ by the commutator of half-twists
along two paths intersecting transversely in one point: $\tilde{B}_n=
B_n/[X_2,X_3^{-1}X_1^{-1}X_2X_1X_3]$. The maps $\sigma:B_n\to S_n$ (induced
permutation) and $\delta:B_n\to \Z$ (degree) factor through
$\tilde{B}_n$, so one can define the subgroups
$\tilde{P}_n=\mathrm{Ker}\,\sigma$ and $\tilde{P}_{n,0}=\mathrm{Ker}\,
(\sigma,\delta)$. The structure of $\tilde{B}_n$ and its subgroups is 
described in detail in \S 1 of \cite{MSegre}; unlike $P_n$ and $P_{n,0}$
which are quite complicated, these groups are fairly easy to understand:
$\tilde{P}_{n,0}$ is solvable, its commutator subgroup is
$[\tilde{P}_{n,0},\tilde{P}_{n,0}]\simeq\Z_2$ and its abelianization
is $\mathrm{Ab}(\tilde{P}_{n,0})\simeq\Z^{n-1}$ (it can in fact be
identified naturally with the reduced regular representation $\mathcal{R}_n$
of $S_n$).
More precisely, we have:

\begin{lemma}[Moishezon]
Let $x_i$ be the image of $X_i$ in $\tilde{B}_n$, and define $s_1=x_1^2$,
$\eta=[x_1^2,x_2^2]$, $u_i=[x_i^{-1},x_{i+1}^2]$ for $1\le i\le n-2$, and 
$u_{n-1}=[x_{n-2}^2,x_{n-1}]$. Then $\tilde{P}_{n,0}$ is
generated by $u_1,\dots,u_{n-1}$, and $\tilde{P}_n$ is generated by
$s_1,u_1,\dots,u_{n-1}$. 

The relations among these elements are
$[u_i,u_j]=1$ if $|i-j|\ge 2$, $[u_i,u_{i+1}]=\eta$, $[s_1,u_i]=1$ if
$i\neq 2$, and $[s_1,u_2]=\eta$. The element $\eta$ is central in
$\tilde{B}_n$, has order $2$ $(i.e.\ \eta^2=1)$, and generates the commutator
subgroups $[\tilde{P}_{n,0},\tilde{P}_{n,0}]=[\tilde{P}_n,\tilde{P}_n]
\simeq \Z_2$ $($in particular, for any two adjacent half-twists $x$ and 
$y$ we have $[x^2,y^2]=\eta)$. 
As a consequence, $\mathrm{Ab}(\tilde{P}_{n})\simeq \Z^n$ and
$\mathrm{Ab}(\tilde{P}_{n,0})\simeq \Z^{n-1}$.

Moreover, the action of $\tilde{B}_n$ on $\tilde{P}_n$ by conjugation is
given by the following formulas: $x_i^{-1}s_1x_i=s_1$ if $i\neq 2$,
$x_2^{-1}s_1x_2=s_1u_2^{-1}$; $x_i^{-1}u_jx_i=u_j$ if $|i-j|\ge 2$,
$x_i^{-1}u_j x_i=u_iu_j$ if $|i-j|=1$, and $x_i^{-1}u_ix_i=u_i^{-1}\eta$.
\end{lemma}

\proof Most of the statement is a mere reformulation of Definition 8 and
Theorem~1 in \S 1.5 of \cite{MSegre}. The only difference is that we define
$u_i$ directly in terms of the generators of $\tilde{B}_n$, while Moishezon 
defines $u_1=(x_2 x_1^2 x_2^{-1}) x_2^{-2}=x_1^{-1}x_2^2 x_1x_2^{-2}$ and
constructs the other $u_i$ by conjugation. In fact, 
$u_i=x^2 y^{-2}$ whenever $x$ and $y$ are
two adjacent half-twists having respectively $i$ and $i+1$
among their end points and such that $xyx^{-1}=x_i$;
our definition of $u_i$ corresponds to the choice $x=x_i^{-1}x_{i+1}x_i$ and 
$y=x_{i+1}$  for $i\le n-2$, and $x=x_{n-2}$
and $y=x_{n-1}x_{n-2}x_{n-1}^{-1}$ for $i=n-1$. Also note that
Moishezon's formula for $x_2^{-1}s_1x_2$ is inconsistent,
due to a mistake in equation (1.25) of \cite{MSegre};
the formula we give is corrected.
\endproof

Intuitively speaking, the reason why $\tilde{B}_n$ is a fairly small group
is that, due to the extra commutation relations, very little is remembered
about the path supporting a given half-twist, namely just its two endpoints 
and the total number of
times that it circles around the $n-2$ other points. This
can be readily checked on simple examples (e.g., half-twists exchanging
the first two points along a path that encircles only one of the $n-2$ other 
points:
since these differ by conjugation by half-twists along paths presenting
a single transverse intersection, they represent the same element in
$\tilde{B}_n$). More generally, we have the following fact:

\begin{lemma}
The elements of $\tilde{B}_n$ corresponding to half-twists exchanging the
first two points are exactly those of the form $x_1 u_1^k \eta^{k(k-1)/2}$
for some integer $k$.
\end{lemma}

\proof Any half-twist exchanging the first two points can be put in the
form $\gamma x_1 \gamma^{-1}$, where $\gamma\in\tilde{P}_n$ can be
expressed as $\gamma=s_1^{\alpha}u_1^{\beta_1}\cdots
u_{n-1}^{\beta_{n-1}}\eta^\epsilon$. Using Lemma 3.1, we have
$x_1^{-1}\gamma x_1=s_1^\alpha
(u_1^{-1}\eta)^{\beta_1}(u_1u_2)^{\beta_2}u_3^{\beta_3}\cdots
u_{n-1}^{\beta_{n-1}}\eta^{\epsilon}$.
Since $(u_1u_2)^{\beta_2}=\eta^{\beta_2(\beta_2-1)/2}u_1^{\beta_2}
u_2^{\beta_2}$, we can rewrite this equality as
$x_1^{-1}\gamma x_1=u_1^{-2\beta_1}\eta^{\beta_1}u_1^{\beta_2}
\eta^{\beta_2(\beta_2-1)/2}\gamma=u_1^k \eta^{k(k-1)/2}\gamma$, where
$k=\beta_2-2\beta_1$. Multiplying by $x_1$ on the left and $\gamma^{-1}$
on the right we obtain $\gamma x_1\gamma^{-1}=x_1 u_1^k\eta^{k(k-1)/2}$.
\endproof

\begin{lemma}
Let $x,y\in\tilde{B}_n$ be elements corresponding to half-twists along
paths with mutually disjoint endpoints. Then $[x,y]=1$.
\end{lemma}

\proof The result is trivial when the paths corresponding to $x$ and $y$
are disjoint or intersect only once. In general, after conjugation we
can assume that $x=\gamma x_1\gamma^{-1}$ for some $\gamma\in\tilde{P}_n$,
and $y=x_3$. By Lemma 3.2, $x=x_1u_1^k\eta^{k(k-1)/2}$ for some
integer $k$. Since $x_1$, $u_1$ and $\eta$ all commute with $x_3$, we
conclude that $[x,y]=1$ as desired.
\endproof

\begin{lemma}
Let $x,y\in\tilde{B}_n$ be elements corresponding to half-twists along
paths with one common endpoint. Then $xyx=yxy$.
\end{lemma}

\proof After conjugation we can assume that $x=x_1$ and $y=\gamma
x_2\gamma^{-1}$ for some $\gamma\in\tilde{P}_n$. By the classification of
half-twists in $\tilde{B}_n$ (Lemma 3.2), there exists an integer $k$ such
that $y=x_2 u_2^k \eta^{k(k-1)/2}=x_2 (s_1u_2^{-1})^{-k} s_1^{k}=s_1^{-k} x_2
s_1^k$. Therefore $xyx=x_1 s_1^{-k} x_2 s_1^k x_1=s_1^{-k} (x_1x_2x_1)
s_1^k=s_1^{-k} (x_2x_1x_2) s_1^k = yxy$.
\endproof

It must be noted that Lemmas 3.3 and 3.4 have also been obtained by Robb
\cite{Robb2}.

\begin{lemma}
The group $\tilde{B}_n$ admits automorphisms $\epsilon_i$ such that
$\epsilon_i(x_i)=x_iu_i$ and $\epsilon_i(x_j)=x_j$ for every $j\neq i$.
Moreover, $\epsilon_i(u_i)=u_i\eta$ and $\epsilon_i(u_j)=u_j$ $\forall j\neq i$.
\end{lemma}

\proof By Lemmas 3.3 and 3.4, 
the half-twists $x_1,\dots,x_{i-1},(x_iu_i),$ $x_{i+1},\dots,x_{n-1}$
satisfy exactly the same relations as the standard generators of
$\tilde{B}_n$. So
$\epsilon_i$ is a well-defined group homomorphism from $\tilde{B}_n$
to itself, and it is injective. The formulas for $\epsilon_i(u_i)$ and 
$\epsilon_i(u_j)$ are easily checked. The surjectivity of $\epsilon_i$
follows from the identity $\epsilon_i(x_iu_i^{-1}\eta)=x_i$.
\endproof

The following definition is motivated by the very particular structure
of the fundamental groups of branch curve complements computed by Moishezon 
for generic projections of $\CP^1\times\CP^1$ and $\CP^2$ 
\cite{MSegre,MVeronese}, which seems to be a feature common to a much larger
class of examples (see \S 4):

\begin{definition}
Define $\tilde{B}_n^{(2)}=\{(x,y)\in \tilde{B}_n\times\tilde{B}_n,
\ \sigma(x)=\sigma(y)\ \mathrm{and}\ \delta(x)=\delta(y)\}$.
We say that the group $\pi_1(\C^2-D_k)$ satisfies
property $(*)$ if there exists an isomorphism 
$\psi$ from $\pi_1(\C^2-D_k)$ to a quotient of
$\tilde{B}_n^{(2)}$ such that, for any geometric generator
$\gamma\in\Gamma_k$, there exist two half-twists $x,y\in \tilde{B}_n$ 
such that $\sigma(x)=\sigma(y)=\theta_k(\gamma)$
and $\psi(\gamma)=(x,y)$. 
\end{definition}

In other words, $\pi_1(\C^2-D_k)$ satisfies property $(*)$ if there exists
a surjective homomorphism from $\tilde{B}_n^{(2)}$ to $\pi_1(\C^2-D_k)$
which maps pairs of half-twists to geometric generators, in a manner
compatible with the $S_n$-valued homomorphisms $\sigma$ and $\theta_k$.

\begin{remark}
If $\pi_1(\C^2-D_k)$ satisfies property $(*)$, then the kernel of
the homomorphism $\theta^+_k:\pi_1(\C^2-D_k)\to S_n\times \Z$ 
is a quotient of $\tilde{P}_{n,0}\times\tilde{P}_{n,0}$ and therefore
a solvable group; in particular
its commutator subgroup is a quotient of $(\Z_2)^2$,
and its abelianization is a quotient of $\,\Z^2\otimes \mathcal{R}_n
\simeq(\Z\oplus\Z)^{n-1}$.
\end{remark}

As an immediate consequence of Definition 3.6
and Lemma 3.3, we have:

\begin{proposition}
If $\pi_1(\C^2-D_k)$ satisfies property $(*)$, then the
stabilization operation is trivial, i.e.\ $K_k=\{1\}$, $G_k=
\pi_1(\C^2-D_k)$, and $G_k^0=\mathrm{Ker}\,\theta^+_k$.
\end{proposition}

\proof Let $\gamma,\gamma'\in\Gamma_k$ be such that $\theta_k(\gamma)$
and $\theta_k(\gamma')$ are disjoint transpositions. Consider the
isomorphism $\psi$ given by Definition 3.6:
there exist half-twists $x,x',y,y'\in\tilde{B}_n$ 
such that $\psi(\gamma)=(x,y)$ and $\psi(\gamma')=(x',y')$.
Since $\theta_k(\gamma)=\sigma(x)=\sigma(y)$ and $\theta_k(\gamma')=
\sigma(x')=\sigma(y')$ are disjoint transpositions, $x$ and $x'$ have
disjoint endpoints, and similarly for $y$ and $y'$. Therefore, by
Lemma 3.3 we have $[x,x']=1$ and $[y,y']=1$, so that $[\psi(\gamma),
\psi(\gamma')]=1$, and therefore $[\gamma,\gamma']=1$. We conclude that
$K_k=\{1\}$, which ends the proof.
\endproof

Let $D_{p,q}$ be the branch curve of a generic polynomial map
$\CP^1\times\CP^1\to\CP^2$ of bidegree $(p,q)$, $p,q\geq 2$.
As will be shown in \S 4, it follows from the computations in 
\cite{MSegre} that $\pi_1(\C^2-D_{p,q})$ satisfies property $(*)$.
This property also holds for the complement of the branch curve of a
generic polynomial map from $\CP^2$ to itself in degree $\geq 3$, 
as follows from the calculations in \cite{MVeronese} 
(see also \cite{TVeronese}), and in various other examples as well 
(see \S 4). It is an interesting question to 
determine whether this remarkable structure of branch curve complements 
extends to generic high-degree projections of arbitrary algebraic surfaces;
this would tie in nicely with a conjecture of Teicher about the virtual 
solvability of these fundamental groups \cite{Tstruct}, and would also
imply Conjecture 1.3.

\section{Examples}

As follows from pp.\ 696--700 of \cite{D1}, if the symplectic manifold $X$
happens to be K\"ahler, then all approximately holomorphic constructions 
can actually be carried out using genuine holomorphic sections of $L^{\otimes
k}$ over $X$, and as a consequence the $\CP^2$-valued maps given by
Theorem 1.1 coincide up to isotopy with projective maps defined by generic
holomorphic sections of $L^{\otimes k}$; therefore, in the case
of complex projective surfaces all calculations can legitimately be 
performed within the framework of complex algebraic geometry.

The fundamental groups of complements of branch curves have already been 
computed for generic projections of various complex projective surfaces.
In many cases, these computations only hold for specific linear systems, 
and do not apply to the high degree situation that we wish to consider.

Nevertheless, it is worth mentioning that, if $D\subset\CP^2$ is the 
branch curve of a generic linear projection of a hypersurface of degree $n$
in $\CP^3$, then it has been shown by Moishezon that 
$\pi_1(\C^2-D)\simeq B_n$ \cite{Mhypersurf}. In fact, in this
specific case there is a well-defined geometric monodromy representation 
morphism $\theta_B$ with values in the braid group $B_n$ rather than 
in the symmetric group $S_n$ as usual, because the $n$ preimages of any
point in $\CP^2-D$ lie in a fiber of the projection 
$\CP^3-\{\mathrm{pt}\}\to\CP^2$, which after trivialization over an
affine subset can be identified with $\C$. Moishezon's computations
then show that $\theta_B:\pi_1(\C^2-D)\to B_n$ is an isomorphism.
An attempt to quotient out $B_n$ by commutators as in the definition of
stabilized fundamental groups yields $\tilde{B}_n$: in this case the
stabilization operation is non-trivial. However this situation is
specific to the linear system $O(1)$, and one expects the fundamental groups
of branch curve complements to behave differently when one instead considers
projections given by sections of $O(k)$ for $k\gg 0$.

Moishezon's result about hypersurfaces in $\CP^3$ has been extended by Robb
to the case of complete intersections (still considering only linear
projections to $\CP^2$ rather than arbitrary linear systems) \cite{Robb2}. 
The result is that, if $D$ is the branch curve for a complete intersection 
of degree $n$ in $\CP^m$ ($m\ge 4$), then the group $\pi_1(\C^2-D)$ is 
isomorphic to $\tilde{B}_n$. It is worth noting that, in this example,
the stabilization operation is trivial. In fact, the groups $\pi_1(\C^2-D)$
can be shown to have property $(*)$ (observe that $\tilde{B}_n$ is the
quotient of $\tilde{B}_n^{(2)}$ by its subgroup $1\times \tilde{P}_{n,0}$).

Conjecture 1.6 holds for $k=1$ in these 
two families of examples: we have $\mathrm{Ab}\,G^0\simeq \Z^{n-1}$ and
$[G^0,G^0]\simeq \Z_2$ in both cases, while $\Z^2/\Lambda_1\simeq \Z$
because the canonical class is proportional to the hyperplane class 
which is primitive.

More interestingly for our purposes, the calculations have also been carried
out in the case of arbitrarily positive linear systems
by Moishezon for two fundamental examples: $\CP^1\times\CP^1$
\cite{MSegre}, and $\CP^2$ \cite{MVeronese} (unpublished, see also 
\cite{TVeronese} for a summary).

\begin{theorem}[Moishezon]
Let $D_{p,q}$ be the branch curve of a generic polynomial map
$\CP^1\times\CP^1\to\CP^2$ of bidegree $(p,q)$, $p,q\geq 2$.
Then the group $\pi_1(\C^2-D_{p,q})$ satisfies property $(*)$,
and its subgroup $H_{p,q}^0=\mathrm{Ker}\,\theta_{p,q}^+$ has the
following structure: $\mathrm{Ab}\,H_{p,q}^0$ is isomorphic to 
$(\Z_2\oplus Z_{p-q})^{n-1}$ if $p$ and $q$ are even, and
$(\Z_{2(p-q)})^{n-1}$ if $p$ or $q$ is odd $($here $n=2pq)$; the commutator subgroup
$[H_{p,q}^0,H_{p,q}^0]$ is isomorphic to $\Z_2\oplus \Z_2$ when $p$ and 
$q$ are even, and $\Z_2$ if $p$ or $q$ is odd.
\end{theorem}

In fact, Moishezon identifies $\pi_1(\C^2-D_{p,q})$ with a quotient of the 
semi-direct product $\tilde{B}_n\ltimes\tilde{P}_{n,0}$,
where $\tilde{B}_n$ acts from the right on $\tilde{P}_{n,0}$ by conjugation
\cite{MSegre}. However it is easy to observe that the map
$\kappa:\tilde{B}_n\ltimes\tilde{P}_{n,0}\to\tilde{B}_n^{(2)}$
defined by $\kappa(x,u)=(x,xu)$ is a group isomorphism (recall the group
structure on $\tilde{B}_n\ltimes\tilde{P}_{n,0}$ is given by
$(x,u)(x',u')=(xx',x'{}^{-1}ux'u')$). The factor $\tilde{P}_{n,0}$ of the
semi-direct product corresponds to the normal subgroup $1\times 
\tilde{P}_{n,0}$ of $\tilde{B}_n^{(2)}$, while the factor $\tilde{B}_n$
corresponds to the diagonally embedded subgroup $\tilde{B}_n=\{(x,x)\}\subset
\tilde{B}_n^{(2)}$.

Moreover, by carefully going over the various formulas identifying a set
of geometric generators for $\pi_1(\C^2-D_{p,q})$ with certain
specific elements in $\tilde{B}_n\ltimes\tilde{P}_{n,0}$
(Propositions 8 and 10 of \cite{MSegre}; cf.\ also \S 1.4, Definition 24
and Remarks 28--29 of \cite{MSegre}), 
or equivalently in $\tilde{B}_n^{(2)}$ after applying the isomorphism
$\kappa$, it is relatively easy to check that each geometric generator 
corresponds to a pair of half-twists with the expected end points
in $\tilde{B}_n^{(2)}$ (see also \S 6 for more details).
Therefore, property $(*)$ and Conjecture 1.3 hold for these groups.

Conjecture 1.6 also holds for $\CP^1\times\CP^1$. Indeed,
$H_2(\CP^1\times\CP^1,\Z)$ is generated by classes $\alpha$ and $\beta$ 
corresponding to the two factors; the hyperplane section class is $L=p\alpha+
q\beta$, while the ramification curve is $R=3L+K=(3p-2)\alpha+(3q-2)\beta$. 
Therefore, the subgroup $\Lambda_{p,q}$ of $\Z^2$ is generated by $(\alpha
\cdot L, \alpha\cdot R)=(q,3q-2)$ and $(\beta\cdot L,\beta\cdot R)=
(p,3p-2)$. An easy computation shows that the quotient $\Z^2/\Lambda_{p,q}=
\Z^2/\langle (q,3q-2),\,
(p,3p-2)\rangle \simeq \Z^2/\langle (q,2),\,(p,2)\rangle$ is isomorphic
to $\Z_2\oplus\Z_{p-q}$ when $p$ and $q$ are even, and to $\Z_{2(p-q)}$
otherwise.

It is worth noting that this nice description for $p,q\ge 2$ completely 
breaks down in the insufficiently ample case $p=1$, where it follows from 
computations of Zariski \cite{Zariski} that $\pi_1(\C^2-D_{1,q})\simeq 
B_{2q}$. So both Conjecture 1.3 and Conjecture 1.6 require a sufficient
amount of ampleness in order to hold ($p,q\ge 2$).

\begin{theorem}[Moishezon]
Let $D_k$ be the branch curve of a generic polynomial map
$\CP^2\to\CP^2$ of degree $k\ge 3$.
Then the group $\pi_1(\C^2-D_k)$ satisfies property $(*)$,
and its subgroup $H_k^0=\mathrm{Ker}\,\theta_k^+$ has the
following structure: $\mathrm{Ab}\,H_k^0$ is isomorphic to 
$(\Z\oplus\Z_3)^{n-1}$ if $k$ is a multiple of $3$, and to $\Z^{n-1}$
otherwise $($here $n=k^2)$; the commutator subgroup
$[H_k^0,H_k^0]$ is trivial for $k$ even and isomorphic to $\Z_2$ for $k$ odd.
\end{theorem}

In this case too, Moishezon in fact identifies $\pi_1(\C^2-D_k)$ with a 
quotient of $\tilde{B}_n\ltimes\tilde{P}_{n,0}$ \cite{MVeronese}
(see also \cite{TVeronese}). Property $(*)$ and Conjecture 1.3 hold for 
$\CP^2$ when $k\ge 3$, but for $k=2$ the group $\pi_1(\C^2-D_2)$ is much 
larger. 

Since $H_2(\CP^2,\Z)$ is generated by the class of a line,
$\Lambda_k$ is the subgroup of $\Z^2$ generated by $(k,3k-3)$, and
$\Z^2/\Lambda_k$ is isomorphic to $\Z\oplus \Z_3$ when $k$ is a multiple
of $3$ and to $\Z$ otherwise. Therefore Conjecture 1.6 holds for $\CP^2$
when $k\ge 3$.

Results for certain projections of Del Pezzo and K3 surfaces have also
been announced by Robb in \cite{Robb2}. 

\begin{theorem}[Robb]
Let $X$ be either a cubic hypersurface in $\CP^3$ or a $(2,2)$ complete
intersection in $\CP^4$, and let $D_k$ be the branch curve of a generic 
algebraic map $X\to\CP^2$ given by sections of $O(kH)$, where $H$ is the
hyperplane section and $k\ge 2$. Then the subgroup 
$H_k^0=\mathrm{Ker}\,\theta_k^+$ of $\pi_1(\C^2-D_k)$ has abelianization 
$\mathrm{Ab}\,H_k^0\simeq \Z^{n-1}$.
\end{theorem}

\begin{theorem}[Robb]
Let $X$ be a K3 surface realized either as a degree $4$ hypersurface in 
$\CP^3$, a $(3,2)$ complete intersection in $\CP^4$ or a $(2,2,2)$ complete
intersection in $\CP^5$, and let $D_k$ be the branch curve of a generic 
algebraic map $X\to\CP^2$ given by sections of $O(kH)$, where $H$ is the
hyperplane section and $k\ge 2$. Then the subgroup 
$H_k^0=\mathrm{Ker}\,\theta_k^+$ of $\pi_1(\C^2-D_k)$ has abelianization 
$\mathrm{Ab}\,H_k^0\simeq (\Z\oplus\Z_k)^{n-1}$.
\end{theorem}

Although to our knowledge no detailed proofs of Theorems 4.3 and 4.4 
have appeared yet, it appears very likely from the sketch of argument given
in \cite{Robb2} that property $(*)$ and Conjecture 1.3 will hold for these
examples as well. In any case we can compare Robb's results with the answers
predicted by Conjecture 1.6.

In the case of the Del Pezzo surfaces, the hyperplane class $H$ is
primitive, and $K=-H$ (so $R_k=(3k-1)H$), so that the subgroup 
$\Lambda_k\subset\Z^2$ is generated by $(k,3k-1)$, and $\Z^2/\Lambda_k
\simeq \Z$, which is in agreement with Theorem 4.3. In the case of the
K3 surfaces, the hyperplane class $H$ is again primitive, but $K=0$
and $R_k=3kH$, so that $\Lambda_k$ is now generated by $(k,3k)$, and
$\Z^2/\Lambda_k\simeq \Z\oplus \Z_k$, in agreement with Theorem 4.4.

The following result for the Hirzebruch surface $\mathbb{F}_1=\mathbb{P}(O_{\CP^1}
\oplus O_{\CP^1}(1))$ is new to our
knowledge; however partial results about this surface have been obtained
by Moishezon, Robb and Teicher \cite{MRTGalHirz,THirzebruch}, and an
ongoing project of Teicher and coworkers is expected to yield another proof
of the same result.

\begin{theorem}
Let $D_{p,q}$ be the branch curve of a generic algebraic map
$\mathbb{F}_1\to\CP^2$ given by three sections of the linear system $O(pF+qE)$,
where $F$ is the class of a fiber, $E$ is the exceptional section, and
$p>q\ge 2$. Then the group $\pi_1(\C^2-D_{p,q})$ satisfies property $(*)$,
and its subgroup $H_{p,q}^0=\mathrm{Ker}\,\theta_{p,q}^+$ has the
following structure: $\mathrm{Ab}\,H_{p,q}^0\simeq (\Z_{3q-2p})^{n-1}$,
where $n=(2p-q)q$, and the commutator subgroup
$[H_{p,q}^0,H_{p,q}^0]$ is isomorphic to 
$\Z_2$ if $p$ is odd and $q$ even, and trivial in all other cases.
\end{theorem}

The proof relies on the observation that $\mathbb{F}_1$ is the blow-up of
$\CP^2$ at one point. Recalling the interpretation of a symplectic (or
K\"ahler) blow-up as the collapsing of an embedded ball, it is easy to
check that $\mathbb{F}_1$ can be degenerated to a union of planes in a manner similar
to $\CP^2$, only with some components missing; most of the calculations 
performed by Moishezon in \cite{MVeronese} for $\CP^2$ can then be re-used
in this context, with the only changes occurring along the exceptional curve
$E$. More details are given in \S 6.2.

As a consequence of property $(*)$, Conjecture 1.3 holds for this example.
So does Conjecture 1.6: indeed, $H_2(\mathbb{F}_1,\Z)$ is generated by $F$ 
and $E$. Recalling that $F\cdot F=0$, $F\cdot E=1$, $E\cdot E=-1$, and
letting $L_{p,q}=pF+qE$ and $R_{p,q}=3L_{p,q}+K=(3p-3)F+(3q-2)E$, we obtain
that $\Lambda_{p,q}\subset\Z^2$ is generated by
$(F\cdot L_{p,q},F\cdot R_{p,q})=(q,3q-2)$ and 
$(E\cdot L_{p,q},E\cdot R_{p,q})=(p-q,3p-3q-1)$. Therefore
$\Z^2/\Lambda_k\simeq \Z^2/\langle (q,3q-2), (p-q,3p-3q-1)\rangle 
\simeq \Z_{3q-2p}$.

A much wider class of examples, including an infinite family of surfaces of
general type, can be investigated if one brings
approximately holomorphic techniques into the picture, although this
makes it only possible to obtain results about the stabilized fundamental
groups of branch curve complements (cf.\ \S 2) rather than the actual 
fundamental groups.

\begin{theorem}
For given integers $a,b\ge 1$ and $p,q\ge 2$,
let $X_{a,b}$ be the double cover of $\CP^1\times\CP^1$ branched along
a smooth algebraic curve of degree $(2a,2b)$, and let
$L_{p,q}$ be the linear system over $X_{a,b}$ defined as the pullback
of $O_{\mathbb{P}^1\times\mathbb{P}^1}(p,q)$ via the double cover.
Let $D_{p,q}$ be the branch curve of a generic approximately holomorphic
perturbation of an algebraic map $X_{a,b}\to\CP^2$ given by three sections 
of $L_{p,q}$. Then the {\em stabilized} fundamental group $G_{p,q}(X_{a,b})=
\pi_1(\C^2-D_{p,q})/K_{p,q}$ satisfies property $(*)$, and its reduced 
subgroup $G^0_{p,q}(X_{a,b})=\mathrm{Ker}\,\theta_{p,q}^+/K_{p,q}$ has
the following structure: $\mathrm{Ab}\,G_{p,q}^0(X_{a,b})\simeq 
(\Z^2/\langle(p,a-2),(q,b-2)\rangle)^{n-1}$,
where $n=4pq$, and the commutator subgroup
$[G_{p,q}^0(X_{a,b}),G_{p,q}^0(X_{a,b})]$ is isomorphic to 
$\Z_2\oplus \Z_2$ if $a,b,p,q$ are all even, trivial if $a$ or $b$ is odd
and $a+p$ or $b+q$ is odd, and isomorphic to $\Z_2$ in all other cases.
\end{theorem}

More precisely, the setup that we consider starts with a holomorphic map 
from $X_{a,b}$ to $\CP^2$ that factors through the double cover 
$X_{a,b}\to\CP^1\times\CP^1$. Such a map is of course not generic in any
sense; however there is a natural explicit way to perturb it in the 
approximately holomorphic category (see \S 7), giving rise to the branch curves 
$D_{p,q}$ that we consider. The map can also be perturbed in the holomorphic
category, which at least for $p$ and $q$ large enough yields a branch curve
that is equivalent to $D_{p,q}$ up to creations and cancellations of pairs 
of nodes. So, on the level of stabilized groups, our
result does give an answer that is relevant from both the symplectic and 
algebraic points of view. Moreover, it is expected that, at least for $p$
and $q$ large enough, the fundamental groups themselves (rather than their
stabilized quotients) should satisfy property $(*)$.

Theorem 4.6 implies that Conjecture 1.6 holds for the manifolds $X_{a,b}$.
Indeed, $X_{a,b}$ can also be described topologically as follows: 
in $\CP^1\times\CP^1$ consider $2a$ curves of the form 
$\CP^1\times\{\mathrm{pt}\}$ and $2b$ curves of the form 
$\{\mathrm{pt}\}\times\CP^1$, and blow up their $4ab$ intersection points
to obtain a manifold $Y_{a,b}$ containing disjoint rational curves
$C_1,\dots,C_{2a}$ (of square $-2b$) and $C'_1,\dots,C'_{2b}$
(of square $-2a$). Then $X_{a,b}$ is the double cover of
$Y_{a,b}$ branched along $C_1\cup\dots\cup C_{2a}\cup C'_1\cup\dots
\cup C'_{2b}$. Now, consider the preimages $\tilde{C}_i=\pi^{-1}(C_i)$
and $\tilde{C}'_i=\pi^{-1}(C'_i)$, and let
$L_{p,q}=p\pi^*\alpha+q\pi^*\beta$ and
$R_{p,q}=3L_{p,q}+K_{X_{a,b}}=(3p+a-2)\pi^*\alpha+(3q+b-2)\pi^*\beta$, 
where $\alpha$ and $\beta$ are the homology generators corresponding to
the two factors of $\CP^1\times\CP^1$.
We have $(\tilde{C}_i\cdot L_{p,q},\tilde{C}_i\cdot R_{p,q})=
(q,3q+b-2)$ and $(\tilde{C}'_i\cdot L_{p,q},\tilde{C}'_i\cdot R_{p,q})=
(p,3p+a-2)$. It is easily shown that these two elements of $\Z^2$
generate the subgroup $\Lambda_{p,q}$; therefore $\Z^2/\Lambda_{p,q}=
\Z^2/\langle (q,3q+b-2),(p,3p+a-2)\rangle\simeq \Z^2/\langle (p,a-2),(q,b-2)
\rangle$.

The techniques involved in the proof of Theorem 4.6, which will be discussed
in \S 7, extend to double covers of other examples for which the answer is
known, possibly including iterated double covers of $\CP^1\times\CP^1$.
One example of particular interest is that of double covers of Hirzebruch
surfaces branched along disconnected curves, for which we make the following
conjecture:

\begin{conj}
Given integers $m,a\ge 1$,
let $X_{2m,a}$ be the double cover of the Hirzebruch surface
$\mathbb{F}_{2m}$ branched along the union of the exceptional
section $\Delta_\infty$ and a smooth algebraic curve 
in the homology class $(2a-1)[\Delta_0]$ $($where $\Delta_0$ is the
zero section, of square $2m)$. Given integers $p,q\ge 2$ such that
$p>2mq$, let $L_{p,q}$ be the linear system over $X_{2m,a}$ defined as 
the pullback of $O_{\mathbb{F}_{2m}}(pF+q\Delta_\infty)$ via the double cover.
Let $D_{p,q}$ be the branch curve of a generic approximately holomorphic
perturbation of an algebraic map $X_{2m,a}\to\CP^2$ given by three sections 
of $L_{p,q}$. Then the reduced stabilized fundamental group 
$G^0_{p,q}(X_{2m,a})=\mathrm{Ker}\,\theta_{p,q}^+/K_{p,q}$ has
abelianization $\mathrm{Ab}\,G_{p,q}^0(X_{2m,a})\simeq 
(\Z^2/\langle(p-2mq,m-2),(2q,2a-4)\rangle)^{n-1}$.
\end{conj}


\section{Stabilized fundamental groups and homological data}

Consider a compact symplectic $4$-manifold $X$ such that $H_1(X,\Z)=0$
and a branched covering map $f_k:X\to\CP^2$ determined by three sections 
of $L^{\otimes k}$, with branch curve $D_k\subset\CP^2$ and geometric 
monodromy representation morphism $\theta_k:\pi_1(\C^2-D_k)\to S_n$. The
purpose of this section is to
construct a natural morphism $\psi_k:\mathrm{Ker}\,\theta_k\to
(\Z^2/\Lambda_k)\otimes\bar\mathcal{R}_n\simeq
(\Z^2/\Lambda_k)^n$ (where $\bar\mathcal{R}_n\simeq \Z^n$ is the regular
representation of $S_n$) and use its properties to prove Theorem 1.5.

Fix a base point $p_0$ in $\C^2-D_k$, and let $p_1,\dots,p_n$ be its
preimages by $f_k$. Let $\gamma\in\pi_1(\C^2-D_k)$ be a loop
in the complement of $D_k$ such that $\theta_k(\gamma)=\mathrm{Id}$. Since 
the monodromy of the branched cover $f_k$ along $\gamma$ is trivial, 
$f_k^{-1}(\gamma)$ is the union of $n$ disjoint closed loops in $X$. 
Denote by $\gamma_i$ the lift of $\gamma$ that starts at the point $p_i$.
Since $H_1(X,\Z)=0$, there exists a surface (or rather a $2$-chain)
$S_i\subset X$ such that $\partial S_i=\gamma_i$. Since $\gamma\subset
\C^2-D_k$, the loop $\gamma_i$ intersects neither the ramification curve
$R_k$ nor the preimage $L_k$ of the line at infinity in $\CP^2$. Therefore,
there exist well-defined algebraic intersection numbers 
$\lambda_i=S_i\cdot L_k$ and $\rho_i=S_i\cdot R_k\in \Z$.
However, there are various possible choices for the surface $S_i$, and
the relative cycle $[S_i]$ is only well-defined up to an element of 
$H_2(X,\Z)$. Therefore, the pair $(\lambda_i,\rho_i)\in\Z^2$ is only
defined up to an element of the subgroup $\Lambda_k$.

\begin{definition}
With the above notations, we denote by $\psi_k:\mathrm{Ker}\,\theta_k
\to (\Z^2/\Lambda_k)^n$ the morphism defined by 
$\psi_k(\gamma)=((S_i\cdot L_k,S_i\cdot R_k))_{1\le i\le n}$.
\end{definition}

In fact, there is no canonical ordering of the preimages of
$p_0$, and $\psi_k$ more naturally takes values in
$(\Z^2/\Lambda_k)\otimes\bar\mathcal{R}_n$, 
as evidenced by Lemma 5.2 below.

Definition 5.1 can naturally be extended to the case 
$H_1(X,\Z)\neq 0$ by instead considering the morphism $\tilde\psi_k:
\mathrm{Ker}\,\theta_k\to H_1(X-L_k-R_k,\Z)^n$
which maps a loop $\gamma$ to the homology classes of its lifts $\gamma_i$
in $X-L_k-R_k$. However, the properties to be expected of this morphism
in general are not entirely clear, due to the lack of available non-simply
connected examples (even though the techniques in \S 6--7 could probably be
applied to the 4-manifold $\Sigma\times\CP^1$ for any Riemann
surface $\Sigma$).

We now investigate the various properties of $\psi_k$.

\begin{lemma}
For every $\gamma\in \mathrm{Ker}\,\theta_k$ and $g\in \pi_1(\C^2-D_k)$,
$\psi_k(g^{-1}\gamma g)=\theta_k(g)\cdot \psi_k(\gamma)$, where $S_n$ acts
on $(\Z^2/\Lambda_k)^n$ by permuting the factors $($i.e., $\psi_k$ is
equivariant$)$.
\end{lemma}

\proof Denoting by $\sigma$ the permutation $\theta_k(g)$, observe that
the lifts of $g^{-1}\gamma g$ are freely homotopic to those of $\gamma$,
and more precisely that the lift of $g^{-1}\gamma g$ through $p_{\sigma(i)}$ 
is freely homotopic to the lift of $\gamma$ through $p_i$. Therefore, the
$\sigma(i)$-th component of $\psi_k(g^{-1}\gamma g)$ is equal to the
$i$-th component of $\psi_k(\gamma)$.
\endproof

\begin{lemma} $K_k\subset \mathrm{Ker}\,\psi_k$, i.e.\ $\psi_k$ factors
through the stabilized group.
\end{lemma}

\proof Recall from Definition 2.2 that $K_k$ is generated by commutators
$[\gamma_1,\gamma_2]$ of geometric generators that are mapped to disjoint
transpositions by $\theta_k$. If $\gamma_1$ is a geometric generator, then
$n-2$ of its lifts to $X$ are contractible closed loops in
$X-L_k-R_k$, while the two other lifts are not closed; and similarly
for $\gamma_2$. However, if $\theta_k(\gamma_1)$ and
$\theta_k(\gamma_2)$ are disjoint, then all the lifts of $[\gamma_1,\gamma_2]$
are contractible loops in $X-L_k-R_k$; therefore $[\gamma_1,\gamma_2]\in
\mathrm{Ker}\,\psi_k$.
\endproof

It is worth noting that, similarly, if $\gamma_1$ and $\gamma_2$ are
geometric generators mapped by $\theta_k$ to adjacent (non-commuting)
transpositions, then $(\gamma_1\gamma_2\gamma_1)(\gamma_2\gamma_1
\gamma_2)^{-1}\in\mathrm{Ker}\,\psi_k$ (only one of the lifts of this
loop is possibly non-trivial, but its algebraic linking numbers with $L_k$ 
and $R_k$ are both equal to zero).

\begin{lemma} For any $\gamma\in \mathrm{Ker}\,\theta_k$, the $n$-tuple
$\psi_k(\gamma)=((\lambda_i,\rho_i))_{1\le i\le n}$ has the property that
$(\sum \lambda_i, \sum \rho_i)\equiv (0,\delta_k(\gamma)) \mod \Lambda_k$.
\end{lemma}

\proof $\gamma\in\pi_1(\C^2-D_k)$ is homotopically trivial in $\C^2$, so
there exists a topological disk $\Delta\subset\C^2$ such that $\partial
\Delta=\gamma$. Now observe that $\partial(f_k^{-1}(\Delta))=\sum \gamma_i$; 
therefore $(\sum \lambda_i,\sum \rho_i)$ is equal (mod $\Lambda_k$) to the 
algebraic intersection numbers of $f_k^{-1}(\Delta)$ with $L_k$ and $R_k$. 
We have $f_k^{-1}(\Delta)\cdot L_k=0$ since $f_k^{-1}(\Delta)\subset 
f_k^{-1}(\C^2)=X-L_k$, and $f_k^{-1}(\Delta)\cdot R_k=\Delta\cdot D_k=
\delta_k(\gamma)$.
\endproof

\begin{lemma} For any geometric generator $\gamma\in\Gamma_k$,
$\psi_k(\gamma^2)=((\lambda_i,\rho_i))_{1\le i\le n}$ is given by
$(\lambda_i,\rho_i)=(0,1)$ if $i$ is one of the two indices exchanged by
the transposition $\theta_k(\gamma)$, and $(\lambda_i,\rho_i)=(0,0)$
otherwise.
\end{lemma}

\proof All lifts of $\gamma^2$ are homotopically trivial, except for
two of them which are freely homotopic to each other and circle once around 
the ramification curve $R_k$.
\endproof

\begin{lemma} There exist two geometric generators $\gamma_1,\gamma_2\in
\Gamma_k$ such that $\theta_k(\gamma_1)=\theta_k(\gamma_2)$ and 
$\psi_k(\gamma_1\gamma_2)=
((-1,0),(1,2),(0,0),\dots,(0,0))$.
\end{lemma}

\proof
Consider a generic line $\ell\subset\CP^2$ intersecting $D_k$ transversely
in $d=\deg D_k$ points, and let $\Sigma=f_k^{-1}(L)$. The restriction
$f_{k|\Sigma}:\Sigma\to\ell=\CP^1$ is a connected simple branched cover of 
degree $n$ with $d$ branch points, with monodromy described by the morphism
$\theta_k\circ i_*:\pi_1(\ell-\{d\ \mathrm{points}\})\to S_n$.
It is a classical fact that the moduli space of all connected simple branched
covers of $\CP^1$ with fixed degree and number of branch points is
connected, i.e.\ up to a suitable reordering of the branch points 
we can assume that the monodromy of $f_{k|\Sigma}$ is described by any given
standard $S_n$-valued morphism. 

So we can find an ordered system of generators $\gamma_1,\dots,\gamma_d$ of
the free group $\pi_1(\ell\cap (\C^2-D_k))$ such that $\theta_k(\gamma_1)=
\theta_k(\gamma_2)=(12)$ and all the other transpositions
$\theta_k(\gamma_i)$ for $i\ge 3$ are elements of $S_{n-1}=\mathrm{Aut}\,
\{2,\dots,n\}$. The loop $\gamma_1\gamma_2$ then belongs to
$\mathrm{Ker}\,\theta_k$, and admits only two non-trivial lifts $g_1$ and
$g_2$ in $\Sigma$, those which start in the first two sheets of the 
branched cover. The loops $g_1$ and $g_2$ bound a topological annulus $A$ which
intersects $R_k$ in two points (projecting to the first two intersection
points of $\ell$ with $D_k$). This annulus separates $\Sigma$ into two
components, a ``large'' component consisting of the sheets numbered from
$2$ to $n$, and a disk $\Delta$ corresponding to the first sheet of the
cover, which does not intersect $R_k$ but contains one of the $n$ preimages
of the intersection point of $\ell$ with the line at infinity in $\CP^2$.
The lift $g_1$ bounds $\Delta$ with reversed orientation; since $\Delta\cdot
R_k=0$ and $\Delta\cdot L_k=1$, the first component of
$\psi_k(\gamma_1\gamma_2)$ is $(-1,0)$. The lift $g_2$ bounds $\Delta\cup
A$; since $A\cdot R_k=2$ and $A\cdot L_k=0$, the second component of
$\psi_k(\gamma_1\gamma_2)$ is $(1,2)$.
\endproof

\proof[Proof of Theorem 1.5]
By Lemma 5.4, $\psi_k$ maps the kernel of $\theta_k^+:\pi_1(\C^2-D_k)\to
S_n\times \Z$ into the subgroup $\Gamma=\{(\lambda_i,\rho_i),\ 
\sum\lambda_i=\sum\rho_i=0\}\simeq(\Z^2/\Lambda_k)\otimes\mathcal{R}_n$
of
$(\Z^2/\Lambda_k)^n$. By Lemma 5.3, $\psi_k$ factors through the quotient
$\mathrm{Ker}\,\theta_k^+/K_k=G_k^0(X,\omega)$, and gives rise to a map
$\phi_k:G_k^0(X,\omega)\to \Gamma\simeq (\Z^2/\Lambda_k)\otimes
\mathcal{R}_n\simeq (\Z^2/\Lambda_k)^{n-1}$. Since
$\Gamma$ is abelian, $[G_k^0,G_k^0]\subset \mathrm{Ker}\,\phi_k$, so
$\phi_k$ factors through the abelianization $\mathrm{Ab}\,G_k^0(X,\omega)$,
as announced in the statement of Theorem 1.5.

We now show that $\phi_k$ is surjective, i.e.\ that $\psi_k$ maps
$\mathrm{Ker}\,\theta_k^+$ {\it onto} $\Gamma$. First, let $\gamma$ and
$\gamma'$ be two geometric generators of $\pi_1(\C^2-D_k)$ corresponding
to adjacent transpositions in $S_n$: then
$\gamma^2{\gamma'}^{-2}\in\mathrm{Ker}\,\theta_k^+$, and Lemma 5.5
implies that $\psi_k(\gamma^2{\gamma'}^{-2})$ has only two non-zero entries,
one equal to $(0,1)$ and the other equal to $(0,-1)$. Recalling from \S 2
that $\theta_k$ is surjective, and using Lemma 5.2, by considering suitable 
conjugates of $\gamma^2{\gamma'}^{-2}$ we can find elements $g_{ij}$ of
$\mathrm{Ker}\,\theta_k^+$ such that $\psi_k(g_{ij})$ has only
two non-zero entries, $(0,1)$ at position $i$ and $(0,-1)$ at position $j$.

Next, consider the geometric generators $\gamma_1,\gamma_2$ given by 
Lemma 5.6: the element $\gamma_1\gamma_2^{-1}$ belongs to
$\mathrm{Ker}\,\theta_k^+$, and $\psi_k(\gamma_1\gamma_2^{-1})=((-1,-1),
(1,1),(0,0),\dots,(0,0))$. Therefore $\psi_k(g_{12}\gamma_1\gamma_2^{-1})=
((-1,0),(1,0),(0,0),\dots,(0,0))$. So, using the surjectivity of $\theta_k$ 
and Lemma 5.2, we can find elements $g'_{ij}$ of $\mathrm{Ker}\,\theta_k^+$
such that $\psi_k(g'_{ij})$ has only two non-zero entries,
$(1,0)$ at position $i$ and $(-1,0)$ at position $j$.
We now conclude that $\psi_k(\mathrm{Ker}\,\theta_k^+)=\Gamma$ by observing
that the $2n-2$ elements $\psi_k(g_{in})$ and $\psi_k(g'_{in})$, $1\le i\le
n-1$, generate $\Gamma$.
\endproof

We finish this section by mentioning two conjectures related to Conjecture
1.6. First of all, we mention that Conjecture 1.6 implies a result about the
fundamental groups of Galois covers associated to branched covers of
$\CP^2$. More precisely, given a complex surface $X$ and a generic
projection $X\to\CP^2$ of degree $n$ with branch curve $D_k$, the associated 
Galois cover $\tilde{X}_k$ is obtained by compactification of the $n$-fold 
fibered product of $X$ with itself above $\CP^2$: the complex surface
$\tilde{X}_k$ 
is a degree $n!$ cover of $\CP^2$ branched along $D_k$. Moishezon and Teicher
have constructed many interesting examples of complex surfaces by this
method, and computed their fundamental groups (see e.g.\ \cite{Tsurv},
\cite{THirzebruch}, \cite{MRTGalHirz}). Given an ordered system of geometric 
generators $\gamma_1,\dots,\gamma_d$ of $\pi_1(\C^2-D_k)$, the fundamental 
group $\pi_1(\tilde{X}_k)$ is known to be isomorphic to the quotient of
$\mathrm{Ker}(\theta:\pi_1(\C^2-D_k)\to S_n)$ by the subgroup generated by
$\gamma_1^2$, $\dots$, $\gamma_d^2$, and $\prod\gamma_i$ (see e.g.\
\cite{THirzebruch}, \S 4). 

By Lemma 5.5, the elements $\gamma_i^2$ and their conjugates map under 
$\psi_k$ to elements of $(\Z^2/\Lambda_k)^n$ with only two non-trivial 
entries $(0,1)$; therefore, assuming Conjecture 1.6, quotienting by all 
squares of geometric generators leads to quotienting the image of $\psi_k$
by $\{(0,\rho_i),\ \sum \rho_i\mathrm{\ is\ even}\}\subset (\Z^2/\Lambda_k)^n$.
Because of Lemma 5.4, and observing that $\delta_k$ takes only even values
on $\mathrm{Ker}\,\theta_k$, we are left with only the first factor in each
summand $\Z^2/\Lambda_k$. Moreover, one easily checks that $\psi_k(\prod
\gamma_i)=((1,0),(1,0),\dots,(1,0))\equiv ((1,0),\dots,(1,0),(1-n,d)) \mod
\Lambda_k$; and by Lemma 5.4, the sum of the first factors is always zero,
so we end up with a group isomorphic to $(\Z_{ks})^{n-2}$, where $ks$ is
the divisibility of $L_k$ in $H_2(X,\Z)$. Moreover, if we also assume that 
property $(*)$ holds in addition to Conjecture 1.6, it can easily be checked
that the commutator subgroup $[G_k^0,G_k^0]$ is contained in the subgroup
generated by the $\gamma_i^2$. Therefore, we have the following conjecture,
satisfied by the examples in \S 4:

\begin{conj}
If $X$ is a simply connected complex surface and $k$ is large
enough, then the fundamental group of the Galois cover $\tilde{X}_k$
associated to a generic projection $f_k:X\to\CP^2$ defined by sections of
$L^{\otimes k}$ is
$\pi_1(\tilde{X}_k)=(\Z_{ks})^{n_k-2}$,
where $ks$ is the divisibility of $L_k$ in $H_2(X,\Z)$ and $n_k=\deg f_k$.
\end{conj}

Also, a careful observation of the examples in \S 4 suggests the following possible
structure for the commutator subgroup $[G_k^0,G_k^0]$, which is worth
mentioning in spite of the rather low amount of supporting evidence:

\begin{conj}
If the symplectic manifold $X$ is simply connected and $k$ is large
enough, then the commutator subgroup $[G_k^0,G_k^0]$ is isomorphic to
$\Gamma_1\times\Gamma_2$, where $\Gamma_1=\Z_2$ if $X$ is spin and
$1$ otherwise, and $\Gamma_2=\Z_2$ if $L_k\equiv K_X\mod 2$ and
$1$ otherwise.
\end{conj}

\section{Moishezon-Teicher techniques for ruled surfaces}

\subsection{Overview of Moishezon-Teicher techniques}
Moishezon and Teicher have developed a general strategy, consisting of
two main steps \cite{Msetup,MSegre,Tsurv}, in order to compute the
group $\pi_1(\C^2-D)$ when $D$ is the branch curve of a generic projection 
to $\CP^2$ of a given projective surface $X\subset\CP^N$.
First, one computes the braid factorization (see \S 2) associated to the 
curve $D$. This calculation involves a degeneration of the surface $X$ to a
singular configuration $X_0$ consisting of a union of planes intersecting
along lines in $\CP^N$, and a careful analysis of the ``regeneration''
process which produces the generic branch curve $D$ out of the singular
configuration \cite{Msetup}. As explained in \S 2, the braid factorization
explicitly provides, via the Zariski-Van Kampen theorem, a (rather
complicated) presentation of the group $\pi_1(\C^2-D)$. In a second step,
one attempts to obtain a simpler description by reorganizing the relations 
in a more orderly fashion and by constructing morphisms between subgroups 
of $\pi_1(\C^2-D)$ and groups related to $\tilde{B}_n$. This process is
carried out in \cite{MSegre} for the case $X\simeq\CP^1\times\CP^1$, and in
subsequent papers for other examples.

\subsubsection{Degenerations and braid monodromy calculations}
The starting point of the calculation is a degeneration of the projective
surface $X\subset\CP^N$ to an arrangement $X_0$ of planes in $\CP^N$ 
intersecting along lines. The degeneration process in the case of manifolds
like $\CP^1\times\CP^1$ and $\CP^2$ is described in detail in \cite{Msetup}.
For example, in the case of $\CP^1\times\CP^1$ embedded by the linear system 
$O(p,q)$, one first degenerates the surface $X$ of degree $2pq$ to a sum of 
$q$ copies of $\CP^1\times\CP^1$ embedded by $O(p,1)$ (each of degree $2p$)
inside $\CP^N$; then each of these surfaces is degenerated into $p$ quadric 
surfaces ($\CP^1\times\CP^1$ embedded by $O(1,1)$); finally, each of the
$pq$ quadric surfaces is degenerated into a union of two planes intersecting
along a line. The resulting arrangement can be represented by the 
diagram in Figure 1.

\begin{figure}[ht]
\begin{center}
\setlength{\unitlength}{7mm}
\begin{picture}(6.5,4.5)(-0.5,-0.15)
\multiput(0,0)(1,0){7}{\line(0,1){4}}
\multiput(0,0)(0,1){5}{\line(1,0){6}}
\put(0,3){\line(1,1){1}}
\put(0,2){\line(1,1){2}}
\put(0,1){\line(1,1){3}}
\multiput(0,0)(1,0){3}{\line(1,1){4}}
\put(5,0){\line(1,1){1}}
\put(4,0){\line(1,1){2}}
\put(3,0){\line(1,1){3}}
\put(0.5,4.15){\makebox(0,0)[cb]{1}}
\put(1.5,4.15){\makebox(0,0)[cb]{2}}
\put(3,4.2){\makebox(0,0)[cb]{\dots}}
\put(5.5,4.15){\makebox(0,0)[cb]{$p$}}
\put(-0.3,0.5){\makebox(0,0)[rc]{1}}
\put(-0.3,1.5){\makebox(0,0)[rc]{2}}
\put(-0.3,2.6){\makebox(0,0)[rc]{$\vdots$}}
\put(-0.3,3.5){\makebox(0,0)[rc]{$q$}}
\put(0.05,0){\makebox(0,-0.05)[lt]{\tiny $00$}}
\put(1.05,0){\makebox(0,-0.05)[lt]{\tiny $10$}}
\put(2.05,0){\makebox(0,-0.05)[lt]{\tiny $20$}}
\put(5.05,0){\makebox(0,-0.05)[lt]{\tiny $p$-1,0}}
\put(0.05,1){\makebox(0,-0.05)[lt]{\tiny $01$}}
\put(1.05,1){\makebox(0,-0.05)[lt]{\tiny $11$}}
\put(2.05,1){\makebox(0,-0.05)[lt]{\tiny $21$}}
\put(6.05,1){\makebox(0,-0.05)[lt]{\tiny $p1$}}
\put(0.05,3){\makebox(0,-0.05)[lt]{\tiny 0,$q$-1}}
\put(1.05,4){\makebox(0,-0.05)[lt]{\tiny $1q$}}
\put(2.05,4){\makebox(0,-0.05)[lt]{\tiny $2q$}}
\put(6.05,4){\makebox(0,-0.05)[lt]{\tiny $pq$}}
\multiput(0,0)(1,0){3}{\circle*{0.08}}
\multiput(0,1)(1,0){3}{\circle*{0.08}}
\multiput(1,4)(1,0){2}{\circle*{0.08}}
\put(0,3){\circle*{0.08}}
\put(5,0){\circle*{0.08}}
\put(6,1){\circle*{0.08}}
\put(6,4){\circle*{0.08}}
\end{picture}
\end{center}
\caption{}
\end{figure}

Each triangle in the diagram represents a plane. Each edge separating 
two triangles represents an intersection line $L_i$ between the corresponding
planes; note that the outer edges of the diagram are not part of the 
configuration.
The branch curve for the projection $X_0\to\CP^2$ is an arrangement of lines
in $\CP^2$ (the projections of the various intersection lines $L_i$); however,
in the regeneration process each of these lines acquires multiplicity
$2$, and the vertices where two or more lines intersect in $X_0$ turn into
certain standard local configurations.

Therefore the braid factorization for $D$ can be computed by looking at the
local contributions of the various vertices in the diagram. Since the
regeneration process turns a local configuration into a branch
curve of degree $2m$, where $m$ is the number of edges meeting at the
given vertex, the local contribution of a vertex is naturally described by
a word in the braid group $B_{2m}$. Moreover, because projecting $X_0$
to $\CP^2$ creates extra intersection points between the projections of
the lines $L_i$ whenever they do not intersect in $X_0$ (i.e.\ when they
do not correspond to edges with a common vertex in the diagram), the
branch curve $D$ contains a number of additional nodes besides the local 
vertex configurations. 

The major difficulty is to arrange the various local
configurations and the additional nodes into a single braid factorization
describing the curve $D$: given a linear projection
$\pi:\CP^2-\{\mathrm{pt}\}\to\CP^1$, one needs to fix a base point in
$\CP^1$ and to choose an ordered system of loops in $\CP^1-\mathrm{crit}\,
\pi_{|D}$ in order to obtain a braid factorization. This choice 
determines in particular how the
local braid monodromy (in $B_{2m}$) for each vertex of the grid is embedded
into the braid monodromy of $D$ (in $B_d$, $d=\deg D$).
A careless setup leads to local embeddings $B_{2m}\hookrightarrow B_d$ that 
may be extremely difficult to determine.

An important observation of Moishezon is that the construction has sufficient
flexibility to allow the images in $\CP^2$ of the various lines and
intersection points to be chosen freely. This makes it possible to use the
following very convenient setup \cite{Msetup}. First choose an ordering of the vertices
in the diagram describing $X_0$; for example, for $\CP^1\times\CP^1$
Moishezon chooses an ordering first by row, then by column, starting from
the lower-left corner of the diagram: 
$00$, $10$, $20$,~$\dots$, $01$, $11$,~$\dots$, $pq$. This determines a
lexicographic ordering of the edges of the diagram: observing that each 
line $L_i$ passes through two vertices $v_i$ and $v'_i$ ($v_i<v'_i$),
the ordering is 
given by $L_i<L_j$ iff either $v'_i<v'_j$, or $v'_i=v'_j$ and $v_i<v_j$.
It is then possible to choose a configuration where the projections of 
the lines $L_i$ are given by equations with real coefficients, with slopes
increasing according to the chosen lexicographic ordering, so that
the intersection of the arrangement of lines in $\CP^2$ with a real slice 
$\mathbb{R}^2$ looks as in Figure 2.

\begin{figure}[ht]
\begin{center}
\setlength{\unitlength}{7mm}
\begin{picture}(6.8,5.2)(-0.5,0)
\put(-0.5,0.75){\line(1,0){6.3}}
\put(-0.5,0){\line(3,1){6.3}}
\put(1,0){\line(2,1){4.8}}
\put(2.5,0){\line(1,1){3.3}}
\put(3.75,0){\line(1,2){2.05}}
\put(4.167,0){\line(1,3){1.63}}
\put(4.7,0){\line(1,5){1.05}}
\put(1.75,0.75){\circle*{0.1}}
\put(2.5,0.75){\circle*{0.1}}
\put(4,1.5){\circle*{0.1}}
\put(5,2.5){\circle*{0.1}}
\put(5.33,3.17){\circle*{0.1}}
\put(5.5,4){\circle*{0.1}}
\put(2,1){\makebox(0,0)[rb]{\small first vertex}}
\put(5.3,4){\makebox(0,0)[rb]{\small last vertex}}
\multiput(2.5,2.3)(0.15,0.1){4}{\circle*{0.01}}
\put(6,0.75){\makebox(0,0)[lc]{$L_1$}}
\put(6,2){\makebox(0,0)[lc]{$L_2$}}
\put(6,5.1){\makebox(0,0)[lc]{$L_{\frac{d}{2}}$}}
\multiput(6.3,3)(0,0.15){4}{\circle*{0.01}}
\end{picture}
\qquad\qquad
\begin{picture}(5,5)(-2.5,-1)
\put(-2.5,2.5){\circle*{0.1}}
\put(-2,2.5){\circle*{0.1}}
\put(-1.25,2.5){\circle*{0.1}}
\put(-0.75,2.5){\circle*{0.1}}
\put(2,2.5){\circle*{0.1}}
\put(2.5,2.5){\circle*{0.1}}
\put(-2.5,2.5){\circle{0.4}}
\put(-2,2.5){\circle{0.4}}
\put(-1.25,2.5){\circle{0.4}}
\put(-0.75,2.5){\circle{0.4}}
\put(2,2.5){\circle{0.4}}
\put(2.5,2.5){\circle{0.4}}
\put(0,0){\line(-1,1){2.36}}
\put(0,0){\line(-4,5){1.87}}
\put(0,0){\line(-1,2){1.16}}
\put(0,0){\line(-1,3){0.77}}
\put(0,0){\line(1,1){2.36}}
\put(0,0){\line(4,5){1.87}}
\multiput(0.5,2.5)(0.15,0){4}{\circle*{0.01}}
\put(-2.5,2.8){\makebox(0,0)[cb]{$\gamma_1$}}
\put(-2,2.8){\makebox(0,0)[cb]{$\gamma'_1$}}
\put(-1.25,2.8){\makebox(0,0)[cb]{$\gamma_2$}}
\put(-0.75,2.8){\makebox(0,0)[cb]{$\gamma'_2$}}
\put(2.55,2.75){\makebox(0,0)[cb]{$\gamma'_{\frac{d}{2}}$}}
\put(2,2.8){\makebox(0,0)[cb]{$\gamma_{\frac{d}{2}}$}}
\end{picture}
\end{center}
\caption{}
\end{figure}

The choice of the slopes of the lines ensures that the intersection points
of $D$ with the reference fiber of $\pi$ (chosen to be $\{x=A\}$ for some 
real number $A\gg 0$) are ordered in the natural way along the real axis,
thus yielding a natural set of geometric generators $\{\gamma_i,\gamma'_i\}$
for $\pi_1(\C^2-D)$, as shown on the right of Figure~2; recall that each line
$L_i$ has multiplicity $2$ and hence yields two generators, and note that the
correct ordering of these generators counterclockwise around the base point
is $\gamma'_{d/2},\gamma_{d/2},\dots,\gamma'_1,\gamma_1$.
Moreover, the various vertices of the diagram describing $X_0$ appear, in 
sequence, for increasing values of $x$ (from left to right).

Since all the contributions to the braid monodromy of $D$ are now localized
along the real $x$-axis, it is a fairly straightforward task to choose a set
of generating loops in the base $\CP^1$ of the fibration $\pi$ and 
enumerate accordingly the various contributions to the braid monodromy
of $D$ (standard configurations at the vertices of the diagram and extra 
nodes coming from pairs of edges without a common vertex). Going through
the list of vertices in decreasing sequence (``from right to left'') yields
the simplest formula (Proposition 1 of \cite{Msetup}):

\begin{proposition}[Moishezon]
With the above setup, the braid monodromy of $D$ is given by the
factorization $\prod_{i=\nu}^1 (C_i\cdot F_i)$, where $\nu$ is the number
of vertices in the diagram, $C_i$ is a product of contributions from
nodal intersections between parts of $D$ corresponding to non-adjacent 
edges, and $F_i$ is the braid monodromy corresponding to the $i$-th vertex,
obtained as the image of a standard local configuration under the embedding
$B_{2m_i}\hookrightarrow B_d$ which maps the standard half-twists generating
$B_{2m_i}$ to half-twists along arcs that remain {\em below} the real
axis.
\end{proposition}

Proposition 6.1 makes it fairly simple to obtain a presentation of 
$\pi_1(\C^2-D)$ in terms of the ``global'' generators 
$\{\gamma_i,\gamma'_i\}$: the nature of the local embeddings 
$B_{2m}\hookrightarrow B_d$ implies that the relations coming from each 
vertex are obtained from standard ``local'' relations (determined by the 
local braid monodromy) simply by renaming each of the $2m$ local geometric 
generators into the corresponding global generator. Additionally, the 
extra nodes yield various commutation relations among geometric generators.

The local configurations for the various types of vertices have been 
analyzed by
Moishezon in \cite{Msetup}, leading to explicit formulas for the local 
contributions to the braid factorization. The easiest case is that of
``2-points'' such as the corner points 00 and $pq$ in the diagram for
$\CP^1\times\CP^1$.
The only line that passes through the vertex locally regenerates to a conic
in $\C^2$, presenting a single vertical tangency near the origin; hence
the local braid monodromy is a single half-twist in $B_2$, giving rise to
an equality relation between the two corresponding geometric generators
of $\pi_1(\C^2-D)$.

The next case is that of ``3-points'' such as those occurring on the
boundary of the diagram for $\CP^1\times\CP^1$. During the first step 
of ``regeneration'', which turns $X_0$ into a union of $pq$ quadric
surfaces, the lines corresponding to the diagonal edges are replaced by 
conics (the branch curve of a bidegree $(1,1)$ map from $\CP^1\times\CP^1$
to $\CP^2$). For the vertices along the top and right sides of the diagram 
(labelled $pj$ or $iq$), the partially regenerated configuration in $\CP^2$
therefore consists of a portion of conic tangent to a line, with the line
having the greatest slope; after further regeneration, the line acquires
multiplicity $2$ and the tangent intersection is replaced by three cusps.
The local contribution to braid monodromy can therefore be expressed by
the product
$\tilde{Z}_{1'2}^3\cdot Z_{1'2'}^3\cdot Z_{1'2}^3\cdot \hat{Z}_{11'}$, where
the various factors are powers of half-twists along the paths represented in
Figure 3 (cf.\ \cite{Msetup} and equation (2.4) in \cite{MSegre}).
The first three factors correspond to cusps arising from the tangent 
intersection between the conic and the line, while the last factor 
corresponds to the vertical tangency of the conic. 
\begin{figure}[ht]
\begin{center}
\setlength{\unitlength}{7mm}
\begin{picture}(4,3)(1,0.6)
\qbezier(5,2)(0,0)(4,3)
\put(1.7,0.637){\line(2,3){1.9}}
\put(6,2){$\longrightarrow$}
\put(5,1.85){\makebox(0,0)[lt]{$1$}}
\put(4.2,3){\makebox(0,0)[lt]{$1'$}}
\put(3.7,3.7){\makebox(0,0)[lt]{$2,2'$}}
\end{picture}
\qquad\qquad\qquad
\begin{picture}(4,3.6)(0,-0.3)
\multiput(0,0)(1,0){4}{\circle*{0.1}}
\multiput(0,1)(1,0){4}{\circle*{0.1}}
\multiput(0,2)(1,0){4}{\circle*{0.1}}
\multiput(0,3)(1,0){4}{\circle*{0.1}}
\qbezier(1,3)(1.75,2.5)(2.5,2.5)
\qbezier(2.5,2.5)(3.5,2.5)(3.5,3)
\qbezier(3.5,3)(3.5,3.3)(3,3.3)
\qbezier(3,3.3)(2.5,3.3)(2,3)
\qbezier(1,2)(2,1.2)(3,2)
\put(1,1){\line(1,0){1}}
\qbezier(0,0)(1,-0.5)(2,-0.5)
\qbezier(2,-0.5)(3.5,-0.5)(3.5,0)
\qbezier(3.5,0)(3.5,0.5)(2.5,0.5)
\qbezier(2.5,0.5)(1.75,0.5)(1,0)
\put(3.8,3){\makebox(0,0)[lc]{$\tilde{Z}_{1'2}^3$}}
\put(3.8,2){\makebox(0,0)[lc]{$Z_{1'2'}^3$}}
\put(3.8,1){\makebox(0,0)[lc]{$Z_{1'2}^3$}}
\put(3.8,0){\makebox(0,0)[lc]{$\hat{Z}_{11'}$}}
\end{picture}
\end{center}
\caption{}
\end{figure}

The 3-points on the bottom and left sides of the diagram give rise to a very
similar local configuration, except for the ordering of the various
components. Finally, the interior vertices of the diagram for 
$\CP^1\times\CP^1$ are all of the same type (``6-points'' in Moishezon's
terminology); a careful analysis of their regeneration yields a certain
braid factorization in $B_{12}$, accounting for the $6$ vertical tangencies,
$24$ nodes and $24$ cusps in the local model, as described in \cite{Msetup}.
The local contributions to the relations defining $\pi_1(\C^2-D)$ have also
been calculated by Moishezon for these various standard models in \S 2 of
\cite{MSegre} (see also below).

\subsubsection{Fundamental group calculations}

The setup described in \S 6.1.1 provides an explicit presentation of
$\pi_1(\C^2-D)$ in terms of geometric generators $\{\gamma_i,\gamma'_i\}$,
$i=1,\dots,\frac{d}{2}$. By Proposition 6.1, the relations consist on one
hand of standard relations given by local models for the various vertices
of the diagram describing the degenerated surface $X_0$, and on the other
hand of commutation relations coming from non-adjacent edges of the diagram.
The goal is then to simplify this presentation and ultimately identify 
$\pi_1(\C^2-D)$ with a certain quotient of $\tilde{B}_n^{(2)}$ (or
$\tilde{B}_n\ltimes\tilde{P}_{n,0}$). In the remainder of this
section, we describe the recipes used by Moishezon for the case
$X=\CP^1\times\CP^1$, following \S 3 of \cite{MSegre}; these methods 
also apply to other complex surfaces admitting similar degenerations, 
such as $X=\CP^2$ \cite{MVeronese} or $X=\mathbb{F}_1$ (\S 6.2).

A first observation of Moishezon is that, after a slight change in the
choice of generators, many of the local relations at the vertices can be
expressed in terms of half of the generators only. More precisely, for each
value of $i$, define a {\it twisting} action $\rho_i$ on the two generators 
$\gamma_i,\gamma'_i$ by the formula $\rho_i(\gamma_i)=\gamma'_i$ and
$\rho_i(\gamma'_i)=\gamma'_i\gamma_i{\gamma'_i}^{-1}$. Choose integers
$l_i$ satisfying the following compatibility conditions: if $i<j$ are the 
labels of the two diagonal edges meeting at a 6-point vertex of the diagram,
then $l_j=l_i-1$; if $i<j$ are the labels of the two vertical edges meeting
at a 6-point, then $l_j=l_i+1$; finally, if $i<j$ are the labels of the two
horizontal edges meeting at a 6-point, then $l_j=l_i$. Now let
$e_i=\rho_i^{l_i}(\gamma_i)$ and $e'_i=\rho_i^{l_i}(\gamma'_i)$.
Because of the {\it invariance properties} of the local models \cite{Msetup},
the local relations corresponding to 2-points and 3-points have the same
expressions in terms of $\{e_i,e'_i\}$ as in terms of
$\{\gamma_i,\gamma'_i\}$, independently of the amount of twisting, 
and those for \mbox{6-points} are also independent of the $l_i$
as long as the compatibility relations hold.
On the other hand, if $i_1<\dots<i_6$ are the labels of the edges meeting 
at a 6-point ($i_1$ and $i_6$ are the two diagonal edges), then it is 
possible to eliminate either $e_{i_1}$ or $e_{i_6}$ from the list of 
generators, because the local
relations imply that
\begin{equation}
e_{i_6}=(e_{i_3}e_{i_2}e_{i_4}^{-1}e_{i_5}^{-1})^{-1}
e_{i_1}(e_{i_3}e_{i_2}e_{i_4}^{-1}e_{i_5}^{-1}).
\end{equation}

The second important observation of Moishezon is that, in many cases
(assuming the diagram is ``large enough'', i.e.\ in the case of a
bidegree $(p,q)$ linear system on $\CP^1\times\CP^1$ that $p,q\ge 2$), 
the relations coming
from cusps and nodes of $D$ can all be reformulated into a very nice
pattern (cf.\ Lemma 14 of \cite{MSegre}). If the two edges $i$ and $j$ 
bound a common triangle in the diagram, then the local relations at their 
common vertex imply that
\begin{equation}
e_ie_je_i=e_je_ie_j,\ e_ie'_je_i=e'_je_ie'_j,\ e'_ie_je'_i=e_je'_ie_j, 
\mbox{ and }e'_ie'_je'_i=e'_je'_ie'_j.
\end{equation}
Otherwise, if there is no
triangle having $i$ and $j$ as edges, or equivalently if the two
transpositions $\theta(e_i)=\theta(e'_i)$ and $\theta(e_j)=\theta(e'_j)\in
S_n$ are disjoint, then we have 
\begin{equation}
[e_i,e_j]=[e_i,e'_j]=[e'_i,e_j]=[e'_i,e'_j]=1.
\end{equation}

Looking at $e_1,\dots,e_{\frac{d}{2}}$, among which there
are only $n-1$ independent generators (by (6.1), many of the
$e_i$ corresponding to diagonal edges can be expressed in terms of the
others), a first consequence of the relations (6.2--6.3) is the following 
(Proposition 8 of \cite{MSegre}):

\begin{lemma}[Moishezon]
In the case of the linear system $O(p,q)$ on $\CP^1\times\CP^1$ $(p,q\ge 2)$,
the subgroup $\mathcal{B}$ of $\pi_1(\C^2-D)$ generated by $e_1,\dots,e_{d/2}$ 
is isomorphic to a quotient of $\smash{\tilde{B}_n}$ $(n=2pq)$.
More precisely, there 
exists a surjective morphism $\tilde\alpha:\tilde{B}_n\to \mathcal{B}$ with 
the property that each $e_i$ is the image of a half-twist in $\tilde{B}_n$, and
$\theta\circ\tilde\alpha=\sigma$ $($i.e.\ the end points of the half-twists
agree with the transpositions $\theta(e_i))$.
\end{lemma}

We now need to add to this description the other generators $e'_i$, or
equivalently the elements $a_i=e'_i e_i^{-1}$.
In the case of $\CP^1\times\CP^1$, we relabel these elements as $d_{ij}$
for the diagonal edge in position $ij$ ($1\le i\le p$, $1\le j\le q$, see
Figure 1), $v_{ij}$ for the vertical edge in position $ij$ ($1\le i<p$,
$1\le j\le q$), and $h_{ij}$ for the horizontal edge in position $ij$ ($1\le
i\le p$, $1\le j<q$). We are especially interested in $a_2=v_{11}$.
Moishezon's next observation is that, as a consequence of relations
(6.2--6.3) and of the local relations of the lower-left-most 6-point in the
diagram, the subgroup generated by $v_{11}$ and the conjugates
$g^{-1}v_{11}g$, $g\in \mathcal{B}$, is naturally isomorphic to a quotient 
of $\tilde{P}_{n,0}$ (\cite{MSegre}, Definition 5 and Lemma 17). Moreover,
the subgroup of $\pi_1(\C^2-D)$ generated by the $e_i$ and by $v_{11}$ is
similarly isomorphic to a quotient of the semi-direct product $\tilde{B}_n
\ltimes\tilde{P}_{n,0}$, or equivalently (as seen in \S 4)
$\tilde{B}_n^{(2)}$. 

The most important relations in $\pi_1(\C^2-D)$ are those coming from the
vertical tangencies of $D$, which we now list for the various types of
vertices. If the edge labelled $i$ passes through a 2-point, then the local
relation $e_i=e'_i$ can be rewritten in the form $a_i=1$. If $i<j$ are
the labels of the two edges meeting at a 3-point, then we have 
$e'_i=e_j^{-1}{e'_j}^{-1}e_ie'_je_j$, or equivalently 
$e'_j=e_i^{-1}{e'_i}^{-1}e_je'_ie_i$. Using (6.2) this relation can be 
rewritten as 
\begin{equation}a_j=e_i^{-1}e_je'_ie_j^{-1}e_ie_j^{-1}=e_i^{-2}(e_ie_j)a_i
(e_j^{-1}e_i^{-1})e_je_i^2e_j^{-1}.\end{equation} Finally, if $i_1<\dots<i_6$
are the labels of the edges meeting at a 6-point (according to the ordering 
rules, $i_1$ and $i_6$ are diagonal, $i_2$ and $i_5$ are vertical, and $i_3$ 
and $i_4$ are horizontal), then, besides (6.1), we also have
\begin{eqnarray}
\begin{cases}
a_{i_6}=(e_{i_3}e_{i_2}e_{i_4}^{-1}e_{i_5}^{-1})^{-1}a_{i_1}
(e_{i_3}e_{i_2}e_{i_4}^{-1}e_{i_5}^{-1})\\
a_{i_5}=(e_{i_1}^{-1}e_{i_3}e_{i_4}^{-1}e_{i_6})^{-1}a_{i_2}
(e_{i_1}^{-1}e_{i_3}e_{i_4}^{-1}e_{i_6})\\
a_{i_4}=(e_{i_1}^{-1}e_{i_2}e_{i_5}^{-1}e_{i_6})^{-1}a_{i_3}
(e_{i_1}^{-1}e_{i_2}e_{i_5}^{-1}e_{i_6})
\end{cases}\\
\begin{cases}
a_{i_3}=(e_{i_3}e_{i_1})^{-1}\,a_{i_2}a_{i_1}(e_{i_1}a_{i_2}^{-1}e_{i_1}^{-1})
\,(e_{i_3}e_{i_1})\\
a_{i_2}=(e_{i_2}e_{i_1})^{-1}\,a_{i_3}a_{i_1}(e_{i_1}a_{i_3}^{-1}e_{i_1}^{-1})
\,(e_{i_2}e_{i_1})
\end{cases}
\end{eqnarray}

A first consequence of relations (6.4--6.6) is that, going inductively
through the various vertices of the grid, all $a_i$ can be expressed in
terms of the $e_1,\dots,e_{d/2}$ and of $a_2=v_{11}$. Therefore 
$\pi_1(\C^2-D)$ is generated by the $e_i$ and by $v_{11}$; hence it is
isomorphic to a quotient of $\tilde{B}_n^{(2)}$. In other words, we
have a surjective homomorphism $\alpha:\tilde{B}_n^{(2)}\to 
\pi_1(\C^2-D)$, extending the morphism $\tilde\alpha:\tilde{B}_n\to
\mathcal{B}$ of Lemma 6.2.%
\medskip

From this point on, the results in \S 3 make it possible to present
Moishezon's argument in a simpler and more illuminating way.
Observe that by Lemma 6.2 each $e_i$ is the image by $\alpha$ of a half-twist
in the diagonally embedded subgroup $\tilde{B}_n\subset \tilde{B}_n^{(2)}$.
Moreover, it is a general fact about irreducible plane curves that all
geometric generators are conjugate to each other in \mbox{$\pi_1(\C^2-D)$};
therefore each of the geometric generators $e_i,e'_i$ is the image of
a pair of half-twists in $\tilde{B}_n^{(2)}$.
Alternately this can be seen directly from the above-listed relations;
these relations also imply that each $a_i$ belongs to the normal subgroup of 
pure degree~$0$ elements $\alpha(\tilde{P}_{n,0}\times
\tilde{P}_{n,0})$, and therefore that the half-twists
corresponding to the geometric generators $e'_i$ have the correct
end points as prescribed by the $S_n$-valued monodromy representation 
morphism $\theta$. Therefore $\pi_1(\C^2-D)$ has the property (*) defined
in \S 3. 

In view of Lemmas 3.3 and 3.4, at this point in the argument
we can discard all the relations in $\pi_1(\C^2-D)$
coming from nodes and cusps of $D$ since they automatically hold
in quotients of $\tilde{B}_n^{(2)}$, and focus on the relations (6.4--6.6)
instead. 

By Lemma 3.2, pairs of half-twists in $\tilde{B}_n^{(2)}$ with
fixed end points can be classified by two integers. More precisely,
fix an ordering of the $n$ sheets of the branched cover $f$, e.g.\ from
left to right and from bottom to top in the diagram. This provides an
ordering of the end points of the half-twists corresponding to $e_i$ and
$e'_i$; we can find an element $g\in\tilde{B}_n^{(2)}$ such that
$e_i=\alpha(g^{-1}(x_1,x_1)g)$,
with ordering of the end points preserved. Then by Lemma 3.2 there exist
integers $k$ and $l$ such that $e'_i=\alpha(g^{-1}(x_1u_1^{-k}\eta^{-k(-k-1)/2},
x_1u_1^{-l}\eta^{-l(-l-1)/2})g)$, i.e.\
$a_i=\alpha(g^{-1}(u_1^k\eta^{k(k-1)/2},u_1^l\eta^{l(l-1)/2})g)$. 
One easily checks by Lemma 3.1 that reversing
the ordering of the end points changes $k$ into $-k$ and $l$ into $-l$.

Since $\alpha$ is a priori not injective,
the integers $k$ and $l$ are not necessarily unique, and there may exist
another pair of integers $(k',l')=(k+\kappa,l+\lambda)$ with the same
property,
i.e.\ such that $\mu=(u_1^\kappa\eta^{k'(k'-1)/2-k(k-1)/2},u_1^\lambda
\eta^{l'(l'-1)/2-l(l-1)/2})\in\mathrm{Ker}\,\alpha$.
If $\kappa$ is odd, then the normal subgroup generated by $\mu$ contains
the commutator of $\mu$ with $(u_2,1)$, which is equal to $(\eta,1)$;
so $(\eta,1)\in\mathrm{Ker}\,\alpha$. If $\kappa$ is even, then 
$\eta^{k'(k'-1)/2-k(k-1)/2}=\eta^{\kappa/2}=\eta^{\kappa(\kappa-1)/2}$
(recall that $\eta^2=1$). Similarly, if $\lambda$ is odd then $(1,\eta)
\in\mathrm{Ker}\,\alpha$, otherwise
$\eta^{l'(l'-1)/2-l(l-1)/2}=\eta^{\lambda(\lambda-1)/2}$. In both cases
we arrive to the conclusion that $\tilde\mu=(u_1^\kappa\eta^{\kappa(\kappa-1)/2},
u_1^\lambda\eta^{\lambda(\lambda-1)/2})\in\mathrm{Ker}\,\alpha$. 
In fact, $\mu$ and $\tilde\mu$ generate the same normal subgroups,
so we also have the converse implication.

Therefore the set of all possible values
for $(\kappa,\lambda)$ forms a subgroup $\Lambda\subset\Z^2$; in fact
$\Lambda=\{(\kappa,\lambda),\ (u_1^\kappa\eta^{\kappa(\kappa-1)/2},
u_1^\lambda\eta^{\lambda(\lambda-1)/2})\in\mathrm{Ker}\,\alpha\}$,
and the pair of integers $(k,l)$ is only defined mod~$\Lambda$.
So, to $e_i$ and $e'_i$ we can associate an element $\bar{a}_i=(k,l)\in
\Z^2/\Lambda$. This element $\bar{a}_i$ contains all the relevant
information about $e_i$ and $e'_i$ apart from the end points. Indeed,
because of Lemma 3.5, up to composition of $\alpha$ with an automorphism
of $\tilde{B}_n^{(2)}$ we can assume $e_i$ to be the image by $\alpha$ of
any given pair of half-twists with the correct end points. And, by Lemma
3.2, if two half-twists $x,y\in \tilde{B}_n$ have the same end points,
then $x^2y^{-2}\in \{1,\eta\}$, so up to a factor of $\eta$ the product
$e'_ie_i=a_ie_i^2$ is determined by $\bar{a}_i$; that ambiguity can
in fact be lifted by arguing that $e_i$ and $e'_i$ are images of 
half-twists.

The subgroup $\Lambda$ can be determined by looking at the relations in
$\pi_1(\C^2-D)$ coming from vertical tangencies of $D$, which determine
the kernel of $\alpha$. We now reformulate these relations in terms of
the $\bar{a}_i$. First, at a 2-point, the relation $a_i=1$ becomes
$\bar{a}_i=(0,0)$. What happens at a 3-point depends on the ordering of 
the sheets of $f$ (i.e., of the triangles of the diagram): the relation
(6.4) becomes
\begin{equation}
\pm \bar{a}_i+\pm \bar{a}_j=(1,1),
\end{equation}
where the first sign is $+$ if the triangle $T$ which has 
both $i$ and $j$ among its edges comes {\it after} the other triangle 
bounded by the edge $i$ and $-$ otherwise, and the second sign is $+$ if 
$T$ comes {\it after} the other triangle bounded by the edge $j$ and $-$ 
otherwise.
In the case of a 6-point with the standard ordering used by Moishezon,
(6.5) and (6.6) become
\begin{equation}
\bar{a}_{i_6}=\bar{a}_{i_1},\qquad
\bar{a}_{i_5}=\bar{a}_{i_2},\qquad
\bar{a}_{i_4}=\bar{a}_{i_3},\qquad
\bar{a}_{i_1}-\bar{a}_{i_2}+\bar{a}_{i_3}=0.
\end{equation}

In the case of $\CP^1\times\CP^1$, denoting by $\bar{d}_{ij}$, $\bar{v}_{ij}$ 
and $\bar{h}_{ij}$ the elements of $\Z^2/\Lambda$ corresponding to $d_{ij}$,
$v_{ij}$ and $h_{ij}$, the relations become (listing the vertices from
left to right and bottom to top): $\bar{d}_{1,1}=(0,0)$,
$\bar{v}_{i,1}-\bar{d}_{i+1,1}=(1,1)$,
$\bar{h}_{1,j}+\bar{d}_{1,j+1}=(1,1)$;
$\bar{d}_{i+1,j+1}=\bar{d}_{i,j}$, $\bar{v}_{i,j+1}=\bar{v}_{i,j}$,
$\bar{h}_{i+1,j}=\bar{h}_{i,j}$,
$\bar{d}_{i,j}-\bar{v}_{i,j}+\bar{h}_{i,j}=0$;
$-\bar{d}_{p,j}-\bar{h}_{p,j}=(1,1)$,
$\bar{d}_{i,q}-\bar{v}_{i,q}=(1,1)$, $\bar{d}_{p,q}=(0,0)$.
Moreover, by construction $\bar{v}_{11}=(0,1)$ (because $v_{11}$ was
identified to a generator of $\tilde{P}_{n,0}$).

Working inductively from the lower-left corner of the diagram, these 
equations yield the formulas 
\begin{equation}\bar{d}_{i,j}=(j-i,0),\qquad
\bar{v}_{i,j}=(1-i,1),\qquad \bar{h}_{i,j}=(1-j,1)\end{equation}
(compare with Proposition
10 of \cite{MSegre}, recalling that the identification between $\tilde{B}_n
\ltimes\tilde{P}_{n,0}$ and $\tilde{B}_n^{(2)}$ is given by
$(x,u)\mapsto (x,xu)$). Moreover, we are left with
the relations $(p-1,-1)=(1,1)$ and $(q-1,-1)=(1,1)$.
In other words, $\Lambda$ is the subgroup of $\Z^2$ generated by $(2-p,2)$
and $(2-q,2)$.

Because all relations in $\pi_1(\C^2-D)$ coming from vertical tangencies
correspond to equality relations between pairs of half-twists in
$\tilde{B}_n^{(2)}$, by the above remarks $\mathrm{Ker}\,\alpha$ is the
normal subgroup of $\tilde{B}_n^{(2)}$ generated by a certain number of
elements of the form $(u_1^\kappa\eta^{\kappa(\kappa-1)/2},u_1^\lambda
\eta^{\lambda(\lambda-1)/2})$, and therefore it is completely determined
by the subgroup $\Lambda\subset\Z^2$.
In our case, $\mathrm{Ker}\,\alpha$ is the normal subgroup of
$\tilde{B}_n^{(2)}$ generated by
$(u_1^{2-p}\eta^{(2-p)(1-p)/2},u_1^2\eta)$ and
$(u_1^{2-q}\eta^{(2-q)(1-q)/2},u_1^2\eta)$. We can now finish the proof
of Theorem 4.1, observing that $H^0_{p,q}=(\tilde{P}_{n,0}\times
\tilde{P}_{n,0})/\mathrm{Ker}\,\alpha$. Recalling from Lemma 3.1 that
$\tilde{P}_{n,0}$ has commutator subgroup $\{1,\eta\}\simeq \Z_2$ and that
$\mathrm{Ab}\,\tilde{P}_{n,0}\simeq \Z^{n-1}$, we have two cases to consider.
First, if e.g.\ $p$ is odd, then by considering the commutator of
$(u_1^{2-p}\eta^{(2-p)(1-p)/2},u_1^2\eta)$ with $(u_2,1)$ we obtain that
$(\eta,1)\in \mathrm{Ker}\,\alpha$ (and similarly if $q$ is odd); but 
one easily checks that $(1,\eta)\not\in\mathrm{Ker}\,\alpha$. On the other
hand, if $p$ and $q$ are both even, then no non-trivial element of
$C=\{1,\eta\}\times\{1,\eta\}$ belongs to $\mathrm{Ker}\,\alpha$. Therefore,
$[H^0_{p,q},H^0_{p,q}]\simeq C/(C\cap \mathrm{Ker}\,\alpha)$ is isomorphic
to $\Z_2$ if $p$ or $q$ is odd, and to $\Z_2\times\Z_2$ if $p$ and $q$ are
even. Moreover, we have $\mathrm{Ab}\,H^0_{p,q}\simeq (\tilde{P}_{n,0}\times
\tilde{P}_{n,0})/\langle C,\mathrm{Ker}\,\alpha\rangle\simeq
(\Z^2/\Lambda)^{n-1}$, which one easily shows to be isomorphic to
$(\Z_2\oplus \Z_{p-q})^{n-1}$ or $(\Z_{2(p-q)})^{n-1}$ depending on the
parity of $p$ and $q$. This completes the proof of Theorem 4.1.
The computations for $\CP^2$ (Theorem 4.2) and other algebraic surfaces 
admitting similar degenerations can be carried out by the same method; for 
example, the case of the Hirzebruch surface $\mathbb{F}_1$ is treated in 
\S 6.2 below.

\begin{figure}[ht]
\begin{center}
\setlength{\unitlength}{7mm}
\begin{picture}(7.5,3.7)(0,-0.7)
\put(0,0){\line(1,1){3}}
\put(0,0){\line(1,0){7}}
\put(1,1){\line(1,0){6}}
\put(2,2){\line(1,0){5}}
\put(3,3){\line(1,0){4}}
\put(7,0){\line(0,1){3}}
\put(4.5,0.5){\makebox(0,0)[cc]{$\mathbb{F}_1$}}
\put(4.5,1.5){\makebox(0,0)[cc]{$\dots$}}
\put(4.5,2.5){\makebox(0,0)[cc]{$\mathbb{F}_1$}}
\put(7.5,1.5){\line(1,0){0.3}}
\put(8,1.5){\line(1,0){0.3}}
\put(8.5,1.47){\makebox(0,0)[lc]{$\to$}}
\end{picture}
\qquad
\begin{picture}(7.5,3.7)(0,-0.7)
\put(1,0){\line(0,1){1}}
\put(2,0){\line(0,1){2}}
\multiput(3,0)(1,0){5}{\line(0,1){3}}
\put(0,0){\line(1,0){7}}
\put(1,1){\line(1,0){6}}
\put(2,2){\line(1,0){5}}
\put(3,3){\line(1,0){4}}
\multiput(0,0)(1,0){5}{\line(1,1){3}}
\put(5,0){\line(1,1){2}}
\put(6,0){\line(1,1){1}}
\put(0.5,-0.7){\makebox(0,0)[cb]{1}}
\put(1.5,-0.7){\makebox(0,0)[cb]{2}}
\put(3.5,-0.7){\makebox(0,0)[cb]{\dots}}
\put(6.5,-0.8){\makebox(0,0)[cb]{$p$}}
\put(7.5,0.5){\makebox(0,0)[lc]{1}}
\put(7.5,1.6){\makebox(0,0)[lc]{$\vdots$}}
\put(7.5,2.5){\makebox(0,0)[lc]{$q$}}
\put(1.05,0){\makebox(0,-0.05)[lt]{\tiny $10$}}
\put(2.05,0){\makebox(0,-0.05)[lt]{\tiny $20$}}
\put(6.05,0){\makebox(0,-0.05)[lt]{\tiny $p$-1,0}}
\put(1.05,1){\makebox(0,-0.05)[lt]{\tiny $11$}}
\put(2.05,1){\makebox(0,-0.05)[lt]{\tiny $21$}}
\put(7.05,1){\makebox(0,-0.05)[lt]{\tiny $p1$}}
\put(3.05,3){\makebox(0,-0.05)[lt]{\tiny $qq$}}
\put(7.05,3){\makebox(0,-0.05)[lt]{\tiny $pq$}}
\put(1,0){\circle*{0.08}}
\put(2,0){\circle*{0.08}}
\put(2,1){\circle*{0.08}}
\put(1,1){\circle*{0.08}}
\put(6,0){\circle*{0.08}}
\put(7,1){\circle*{0.08}}
\put(3,3){\circle*{0.08}}
\put(7,3){\circle*{0.08}}
\end{picture}
\end{center}
\caption{}
\end{figure}

\subsection{The Hirzebruch surface $\mathbb{F}_1$}
In this section, we prove Theorem 4.5 using the method outlined in the
preceding section. Consider the projective embedding of $\mathbb{F}_1$
defined by sections of the linear system $O(pF+qE)$, $p>q\ge 2$ (recall
$F$ is the fiber and $E$ is the exceptional section). This projective
surface can be degenerated in the same manner as the Veronese surface
of which it is a blow-up (the projective embedding of $\CP^2$ defined by
sections of $O(p)$), following the procedure described in \S 3 of
\cite{Msetup}. This surface of degree $n=(2p-q)q$ can be first degenerated
into a sum of $q$ Hirzebruch surfaces, of degrees respectively $2p-1$,
$2p-3$, \dots, $2(p-q)+1$. Each of these Hirzebruch surfaces can then be
degenerated into the union of a plane and a certain number of quadric
surfaces, which in turn can each be degenerated to two planes. The resulting
diagram is pictured in the right half of Figure 4.

One uses the same setup as in \S 6.1.1, ordering the vertices from left to
right and bottom to top, and the edges accordingly. The braid monodromy is
given by Proposition 6.1. It follows from Moishezon's work that all vertices
correspond to well-known configurations: the two vertices $qq$ and $pq$ are 
$2$-points, while the other boundary vertices are
$3$-points and the interior vertices are $6$-points.

As in \S 6.1.2, one replaces the natural set of geometric generators
$\{\gamma_i,\gamma'_i\}$ by twisted generators $e_i=\rho_i^{l_i}(\gamma_i)$
and $e'_i=\rho_i^{l_i}(\gamma'_i)$, where the integers $l_i$ satisfy the
required compatibility conditions, in order to have (6.1) at
all 6-points. Moreover, relations (6.2) and (6.3) hold for all pairs of
edges ((6.2) if the edges bound a common triangle, (6.3) otherwise), by 
the same argument as for $\CP^2$: the proof of Lemma 1 of \cite{MVeronese} 
(see also Lemma 14 of \cite{MSegre}) applies almost without modification.

Eliminating redundant diagonal edges as allowed by (6.1), we are left
with exactly $n-1$ independent generators among the $e_i$. As in the case
of $\CP^1\times\CP^1$, relations (6.2) and (6.3) imply that the subgroup
$\mathcal{B}$ generated by the $e_i$ is isomorphic to a quotient of
$\tilde{B}_n$, and Lemma 6.2 extends to the case of the Hirzebruch surface
$\mathbb{F}_1$.

As previously, we let $a_i=e'_i e_i^{-1}$, and we relabel these elements as
$d_{ij}$, $v_{ij}$ and $h_{ij}$. We are now interested in $a_1=v_{11}$~:
one can again show that the subgroup generated by $v_{11}$ and the
conjugates $g^{-1}v_{11}g$, $g\in\mathcal{B}$ is isomorphic to a quotient of
$\tilde{P}_{n,0}$, by Lemma 5 of \cite{MVeronese} (the argument is the same
for $\mathbb{F}_1$ as for $\CP^2$); the subgroup of $\pi_1(\C^2-D)$
generated by the $e_i$ and by $a_1$ is again isomorphic to a quotient of
$\tilde{B}_n\ltimes\tilde{P}_{n,0}\simeq \tilde{B}_n^{(2)}$. 

Relations (6.4--6.6) imply that, going through the various $3$-points and
$6$-points of the diagram, all the $a_i$ can be expressed in terms of 
$e_1,\dots,e_{d/2}$ and $a_1=v_{11}$; therefore $\pi_1(\C^2-D)$ is generated
by $e_1,\dots,e_{d/2}$ and $a_1$, so that we again obtain a surjective
morphism $\alpha:\tilde{B}_n^{(2)}\to\pi_1(\C^2-D)$. As in the case of
$\CP^1\times\CP^1$, the various geometric generators are images by $\alpha$
of pairs of half-twists with correct end points, so that property (*) holds
once more. Using the classification of half-twists in $\tilde{B}_n$
(Lemma 3.2), we can consider pairs of integers $\bar{a}_i$ instead of the
elements $a_i$; once again, the $\bar{a}_i$ are only defined modulo a
certain subgroup $\Lambda\subset\Z^2$.

The various relations between the $\bar{a}_i$ are now the following:
$\bar{v}_{i,1}-\bar{d}_{i+1,1}=(1,1)$,
$\bar{v}_{i,i}-\bar{h}_{i+1,i}=(1,1)$;
$\bar{d}_{i+1,j+1}=\bar{d}_{i,j}$, $\bar{v}_{i,j+1}=\bar{v}_{i,j}$,
$\bar{h}_{i+1,j}=\bar{h}_{i,j}$,
$\bar{d}_{i,j}-\bar{v}_{i,j}+\bar{h}_{i,j}=0$;
$-\bar{d}_{p,j}-\bar{h}_{p,j}=(1,1)$,
$\bar{v}_{q,q}=(0,0)$,
$\bar{d}_{i,q}-\bar{v}_{i,q}=(1,1)$, $\bar{d}_{p,q}=(0,0)$.
Moreover, $\bar{v}_{1,1}=(0,1)$.
Therefore, $\bar{d}_{i,j}=(2j-2i+1,j-i+1)$,
$\bar{v}_{i,j}=(2-2i,2-i)$ and $\bar{h}_{i,j}=(1-2j,1-j)$
(compare with Proposition 4 of \cite{MVeronese}), and we
are left with two additional relations: $(2p-2,p-2)=(1,1)$ and
$(2-2q,2-q)=(0,0)$. Therefore, $\Lambda$ is the subgroup of $\Z^2$
generated by $(2p-3,p-3)$ and $(2q-2,q-2)$, and $\mathrm{Ker}\,\alpha$
is the normal subgroup of $\tilde{B}_n^{(2)}$ generated by
$(u_1^{2p-3}\eta^{(2p-3)(2p-4)/2},u_1^{p-3}\eta^{(p-3)(p-4)/2})$ and
$(u_1^{2q-2}\eta^{(2q-2)(2q-3)/2},u_1^{q-2}\eta^{(q-1)(q-2)/2})$.

Considering the commutator of the first generator with $(u_2,1)$, we
obtain that $(\eta,1)\in\mathrm{Ker}\,\alpha$. Moreover, if either $p$ is
even or $q$ is odd, then considering the commutator of one of the generators
with $(1,u_2)$, we obtain that $(1,\eta)\in\mathrm{Ker}\,\alpha$. On the
contrary, if $p$ is odd and $q$ is even then
$(1,\eta)\not\in\mathrm{Ker}\,\alpha$. We conclude that
$[H^0_{p,q},H^0_{p,q}]\simeq C/(C\cap \mathrm{Ker}\,\alpha)$ is trivial or 
isomorphic to $\Z_2$ depending on the parity of $p$ and $q$, and that
$\mathrm{Ab}\,H^0_{p,q}\simeq (\Z^2/\Lambda)^{n-1}\simeq (\Z^2/\langle
(p,3),(q,2)\rangle)^{n-1}\simeq (\Z_{3q-2p})^{n-1}$.

\section{Double covers of $\CP^1\times\CP^1$}

In this section, we sketch the proof of Theorem 4.6, which combines
the methods described in \S 6 with ideas similar to those in \cite{AK2}.

\subsection{Generic perturbations of iterated branched covers}
Let $C$ be a smooth algebraic curve of degree $(2a,2b)$ in
$Y=\CP^1\times\CP^1$, and let $X_{a,b}$ be the double cover of $Y$ branched 
along $C$. Then one can construct a map $f^0:X_{a,b}\to\CP^2$ simply by 
composing the double cover $\pi:X_{a,b}\to Y$ with a generic projective map 
$g:Y\to\CP^2$ determined by sections of $O(p,q)$. The map $f^0$ is not generic~:
its ramification curve is the union of the ramification curve of
$\pi$ and the preimage by $\pi$ of the ramification
curve of $g$, and so the branch curve $D^0$ of $f^0$ is the 
union of $g(C)$ (with multiplicity $1$) and the branch curve $D_g$ of $g$ 
(with multiplicity $2$). 

This situation is extremely similar to that considered in
\cite{AK2} for the composition of a generic map from a symplectic 4-manifold
to $\CP^2$ with a quadratic map from $\CP^2$ to itself. The local behavior
of the map $f^0$ is generic everywhere except at the intersection points
of $C$ with the ramification curve of $g$~; assuming that $C$ and $g$ are 
chosen generically, a local model for $f^0$ near these points is 
$(x,y)\mapsto (-x^2+y,-y^2)$, for which a generic local perturbation is given
e.g.\ by $(x,y)\mapsto (-x^2+y,-y^2+\epsilon x)$ where $\epsilon$
is a small non-zero constant (cf.\ also \cite{AK2}). There are several
ways in which the map $f^0$ can be perturbed and made generic. If
the linear system $\pi^*O(p,q)$ is sufficiently ample, then $f^0$ can be
deformed within the holomorphic category into a generic projective map 
which no longer factors through the 
double cover $\pi$. Another possibility, if $p$ and $q$ are sufficiently
large, is to use approximately holomorphic methods (Theorem 1.1) to deform
$f^0$ into a map with generic local models (cf.\ \cite{AK2}). 

In both cases, the effect of the perturbation on the topology of the branch
curve of $f^0$ is pretty much the same. First, the local model near an
intersection point of $C$ with the ramification curve of $g$ is perturbed as
described above (up to isotopy), which transforms a tangent intersection of
$g(C)$ with the branch curve of $g$ in $\CP^2$ into a standard configuration
with three cusps \cite{AK2}. Secondly, the two copies of the branch curve
of $g$, which make up the multiplicity two component of $D^0$, are separated
and made transverse to each other; this deformation of $D_g$
is performed either within the holomorphic category or resorting to
approximately holomorphic perturbations. In the second case, the
perturbation process can be performed in a very flexible manner, which in
some cases may create negative intersections~; restricting oneself to
algebraic perturbations is a convenient way to avoid this phenomenon, but
makes the global perturbation harder to describe explicitly. In any case,
up to isotopy and creation or cancellation of pairs of intersections between
the two deformed copies of the branch curve of $g$, the topology of the 
resulting generic branch curve $D$ is uniquely determined and can be computed 
easily from that of $D^0$. In fact, the approximately holomorphic
perturbation process can always be carried out, even for small values of $p$
and $q$ for which neither the holomorphic construction nor Theorem 1.1 are
able to yield generic projective maps~; in this situation, we can still
study the topology of the curve $D$, but Theorem 4.6 only describes a
``virtual'' generic projective map.

As in \S 6, the study of the curve $D$ relies on a degeneration process:
one first degenerates the curve $C$ in $Y=\CP^1\times\CP^1$ into a union of
two sets of parallel lines, $2a$ along one factor and $2b$ along the other
factor. Parallel lines are then merged, so that the resulting configuration
$C_0\subset Y$ consists of only two
components, a $(1,0)$-line of multiplicity $2a$ and a $(0,1)$-line of
multiplicity $2b$. Finally, one degenerates the projective embedding of $Y$
given by the linear system $O(p,q)$ into 
an arrangement $Y_0$ of planes intersecting along lines, as in \S 6.1.
The fully degenerated branch curve is a union of
lines, some of which correspond to the intersections between the planes in
$Y_0$ (each contributing with multiplicity $4$, since the branch curve of
$g$ is counted with multiplicity $2$), while the others are the images of
the $p+q$ components into which $C_0$ degenerates (some of these components
contribute with multiplicity $2a$, others with multiplicity $2b$).

The curve $D$ can be recovered from this arrangement of lines by the
converse ``regeneration'' process, which first yields the union $D_g\cup
g(C_0)$ (by deforming $Y_0$ into the smooth surface $Y$), then $D_g\cup
g(C)=D^0$ (by separating the multiple components of $C_0$ and smoothing the
resulting curve), and finally $D$ (by performing the prescribed local
perturbation at the intersection points of the two ramification curves
and by perturbing the two copies of $D_g$ in a generic way).

\subsection{Braid monodromy calculations} The braid monodromy for the
curve $D_g\cup g(C_0)$ (and for the subsequent regenerations $D^0$ and $D$)
can be computed using the same methods as in \S 6.1.1. The diagram
describing the degenerated configuration is as represented on Figure 5,
which differs from Figure 1 only by the addition of edges corresponding to 
$C_0$ along the top and right boundaries of the diagram.

\begin{figure}[ht]
\begin{center}
\setlength{\unitlength}{7mm}
\begin{picture}(6.5,4.5)(-0.5,-0.15)
\put(0,4.02){\line(1,0){6}}
\put(0,3.98){\line(1,0){6}}
\put(5.98,0){\line(0,1){4}}
\put(6.02,0){\line(0,1){4}}
\multiput(0,0)(1,0){7}{\line(0,1){4}}
\multiput(0,0)(0,1){5}{\line(1,0){6}}
\put(0,3){\line(1,1){1}}
\put(0,2){\line(1,1){2}}
\put(0,1){\line(1,1){3}}
\multiput(0,0)(1,0){3}{\line(1,1){4}}
\put(5,0){\line(1,1){1}}
\put(4,0){\line(1,1){2}}
\put(3,0){\line(1,1){3}}
\put(0.5,4.15){\makebox(0,0)[cb]{1}}
\put(1.5,4.15){\makebox(0,0)[cb]{2}}
\put(3,4.2){\makebox(0,0)[cb]{\dots}}
\put(5.5,4.15){\makebox(0,0)[cb]{$p$}}
\put(-0.3,0.5){\makebox(0,0)[rc]{1}}
\put(-0.3,1.5){\makebox(0,0)[rc]{2}}
\put(-0.3,2.6){\makebox(0,0)[rc]{$\vdots$}}
\put(-0.3,3.5){\makebox(0,0)[rc]{$q$}}
\put(0.05,0){\makebox(0,-0.05)[lt]{\tiny $00$}}
\put(1.05,0){\makebox(0,-0.05)[lt]{\tiny $10$}}
\put(2.05,0){\makebox(0,-0.05)[lt]{\tiny $20$}}
\put(6.05,0){\makebox(0,-0.05)[lt]{\tiny $p0$}}
\put(0.05,1){\makebox(0,-0.05)[lt]{\tiny $01$}}
\put(1.05,1){\makebox(0,-0.05)[lt]{\tiny $11$}}
\put(2.05,1){\makebox(0,-0.05)[lt]{\tiny $21$}}
\put(6.05,1){\makebox(0,-0.05)[lt]{\tiny $p1$}}
\put(0.05,4){\makebox(0,-0.05)[lt]{\tiny $0q$}}
\put(1.05,4){\makebox(0,-0.05)[lt]{\tiny $1q$}}
\put(2.05,4){\makebox(0,-0.05)[lt]{\tiny $2q$}}
\put(6.05,4){\makebox(0,-0.05)[lt]{\tiny $pq$}}
\multiput(0,0)(1,0){3}{\circle*{0.08}}
\multiput(0,1)(1,0){3}{\circle*{0.08}}
\multiput(0,4)(1,0){3}{\circle*{0.08}}
\put(6,0){\circle*{0.08}}
\put(6,1){\circle*{0.08}}
\put(6,4){\circle*{0.08}}
\end{picture}
\end{center}
\caption{}
\end{figure}

Thanks to Proposition 6.1, we only need to understand the local behavior
of the curves $D_g\cup g(C_0)$, $D^0$ and $D$ near the various vertices
of the diagram. At all vertices except those through which $C_0$ passes
(top and right sides of the diagram), the local description of $D_g\cup
g(C_0)$ and $D^0$ is exactly the same as that of $D_g$, which has already
been discussed in \S 6.1~: the various vertices are standard 2-points, 
3-points and 6-points as in Moishezon's work \cite{MSegre}. 
\begin{figure}[ht]
\begin{center}
\setlength{\unitlength}{7mm}
\begin{picture}(15,4.6)(-3,-0.1)
\put(-0.3,-0.5){\line(1,2){1.15}} 
\put(0.2,-0.55){\line(1,3){0.78}}
\put(1.15,2.2){\line(1,2){1}}
\put(1.02,2.2){\line(1,3){0.67}}
\put(0,1.45){\line(2,1){2.6}}
\put(0,1.55){\line(2,1){2.6}}
\qbezier(2.6,2.75)(3.5,3.25)(2.6,2.85)
\qbezier(3.4,3.15)(2.5,2.75)(3.4,3.25)
\put(3.4,3.15){\line(2,1){2.5}}
\put(3.4,3.25){\line(2,1){2.5}}
\qbezier(0.845,1.8)(1.06,2.1)(0.97,1.79)
\qbezier(1.14,2.2)(0.94,1.9)(1.015,2.2)
\qbezier(0.2,0.5)(0.9,2)(3,2.95)
\qbezier(0.2,0.5)(0,-0.5)(2,-0.5)
\qbezier(3,2.95)(4,3.45)(6,4)
\qbezier(2,-0.5)(4,-0.5)(6,0)
\put(6.1,0){\makebox(0,0)[lc]{$1$}}
\put(6.1,4){\makebox(0,0)[lc]{$1'$}}
\put(6.1,4.5){\makebox(0,0)[lc]{$2,2'$}}
\put(2.25,4.3){\makebox(0,0)[lc]{$3$}}
\put(1.75,4.5){\makebox(0,0)[rc]{$4$}}
\put(1,2){\circle*{0.1}}
\put(0.8,2.2){\makebox(0,0)[rc]{B}}
\put(0.2,0.5){\circle*{0.1}}
\put(0,0.5){\makebox(0,0)[rc]{C}}
\put(3,3){\circle*{0.1}}
\put(3,3.2){\makebox(0,0)[cb]{A}}
\put(4.2,2.8){\makebox(0,0)[cc]{$D_g$}}
\put(1.1,3.8){\makebox(0,0)[rt]{$g(C_0)$}}
\put(-3,4){\line(1,0){2}}
\put(-3,4.02){\line(1,0){2}}
\put(-3,3.98){\line(1,0){2}}
\put(-2,3){\line(0,1){1}}
\put(-3,3){\line(1,1){1}}
\put(-2.8,4.2){\makebox(0,0)[cb]{\small $3$}}
\put(-1.2,4.2){\makebox(0,0)[cb]{\small $4$}}
\put(-2.8,3.25){\makebox(0,0)[rb]{\small $1$}}
\put(-1.9,3.2){\makebox(0,0)[lb]{\small $2$}}
\put(8.4,4){\makebox(0,0)[lb]{B}}
\put(8.5,2.65){\line(3,1){3}}
\put(8.5,2.75){\line(3,1){3}}
\put(11.5,4.2){\makebox(0,0)[lc]{$D^0$}}
\qbezier[350](9.9,4.4)(8.9,1.8)(10.9,4.6)
\qbezier[350](9.6,4.3)(8.6,1.7)(10.6,4.5)
\qbezier[350](10.2,4.5)(9.2,1.9)(11.2,4.7)
\qbezier[350](10.5,4.6)(9.5,2)(11.5,4.8)
\qbezier[350](9,1.65)(10,4.25)(8,1.45)
\qbezier[350](9.3,1.75)(10.3,4.35)(8.3,1.55)
\qbezier[350](9.6,1.85)(10.6,4.45)(8.6,1.65)
\qbezier[350](9.9,1.95)(10.9,4.55)(8.9,1.75)
\put(8.4,0){\makebox(0,0)[lb]{C}}
\put(10.8,0.6){\makebox(0,0)[lc]{$D^0$}}
\qbezier(9.5,-0.5)(10,1)(11.5,1.5)
\put(9,-0.49){\line(1,1){2.5}}
\put(9,0){\line(3,2){2.5}}
\put(9,-0.2){\line(5,4){2.5}}
\put(11.2,1.98){\line(-4,-5){2}}
\end{picture}
\end{center}
\caption{}
\end{figure}
Moreover, the
local configuration for $D$ at such a vertex simply consists of two copies 
of the local configuration for $D_g$, shifted apart from each other by 
a generic translation. The two components, which correspond to the two 
preimages of the ramification curve of $g$ under the branched cover $\pi$,
may intersect at nodal points of either orientation~; we won't be overly 
concerned by the details of these intersections, since the various
possible configurations only differ by isotopies and creations or 
cancellations of pairs of nodes, which do not affect the stabilized
fundamental group in any way.

We now consider a vertex along the top boundary of the diagram, at position
$iq$ with $1\le i\le p-1$. The local configuration for $D_g\cup g(C_0)$ at 
such a point is as shown on Figure 6. The parts labelled $1,1',2,2'$
correspond to $D_g$, and form a standard 3-point (cf.\ \S 6.1.1 and Figure 3),
presenting three cusp singularities near the point A. The parts labelled
$3$ and $4$ correspond to $g(C_0)$, obtained by ``regeneration'' of the two
lines associated to the horizontal edges of the diagram passing through the 
vertex. The curve $g(C_0)$ presents tangent intersections with the two lines 
$2$ and $2'$ near the point B, and with the conic $1,1'$ at the point C.
The two intersections of the line labelled $4$ with the conic $1,1'$ in
$\CP^2$ remain as nodes since the corresponding curves fail to intersect
in $Y$.

The local description of the curve $D^0=D_g\cup g(C)$ is obtained from that
of $D_g\cup g(C_0)$ by separating $C_0$ into $2b$ parallel components~; this
yields $2b$ copies of the lines labelled $3$ and $4$ in Figure 6, and the
local configuration near the points B and C becomes as shown in the right
half of Figure 6 (the pictures correspond to the case $b=2$). Finally, in
order to obtain $D$ we must perturb $D^0$ in the manner explained in \S 7.1:
the multiplicity two component $D_g\subset D^0$ (corresponding to the parts
labelled $1,1',2,2'$ in Figure 6) is separated into two distinct copies
(in particular the point A is duplicated), while each tangent intersection 
of $g(C)$ with $D_g$ (such as those near points B and C) gives rise to three
cusps. It is then possible to write explicitly the local braid monodromy for
$D$, with values in $B_{4b+8}$ by enumerating carefully the $4b+2$ vertical
tangencies, $18b+6$ cusps, and nodes of the local model (the exact number of 
nodes depends on the choice of boundary values for the local perturbation of
$D^0$).

In fact, since we only aim to compute {\it stabilized} fundamental groups
of branch curve complements, we shall not concern ourselves with
the nodes of $D$, since these only yield commutation relations which by
definition always hold in the stabilized group. 
\begin{figure}[ht]
\begin{center}
\setlength{\unitlength}{5mm}
\begin{picture}(22,3.2)(0,0.5)
\multiput(0,3)(1,0){8}{\circle*{0.15}}
\put(0,3.2){\makebox(0,0)[cb]{\small $1$}}
\put(1,3.2){\makebox(0,0)[cb]{\small $\tilde{1}$}}
\put(2,3.2){\makebox(0,0)[cb]{\small $1'$}}
\put(3,3.2){\makebox(0,0)[cb]{\small $\tilde{1}'$}}
\put(4,3.25){\makebox(0,0)[cb]{\small $2$}}
\put(5,3.2){\makebox(0,0)[cb]{\small $\tilde{2}$}}
\put(6,3.2){\makebox(0,0)[cb]{\small $2'$}}
\put(7,3.2){\makebox(0,0)[cb]{\small $\tilde{2}'$}}
\put(7.8,3.1){\makebox(0,0)[lc]{$\kappa^A_{1/2/3}$}}
\qbezier[100](2,3)(3,2.3)(4,3)
\qbezier[70](2,3)(2.7,2.5)(3.4,2.5)
\qbezier[70](6,3)(5.3,2.5)(4.6,2.5)
\put(3.4,2.5){\line(1,0){1.2}}
\qbezier[80](2,3)(2.75,2.5)(3.5,3)
\qbezier[80](6,3)(5.25,2.5)(4.5,3)
\qbezier[80](3.5,3)(4,3.33)(4.5,3)
\multiput(11,3)(1,0){8}{\circle*{0.15}}
\put(11,3.2){\makebox(0,0)[cb]{\small $1$}}
\put(12,3.2){\makebox(0,0)[cb]{\small $\tilde{1}$}}
\put(13,3.2){\makebox(0,0)[cb]{\small $1'$}}
\put(14,3.2){\makebox(0,0)[cb]{\small $\tilde{1}'$}}
\put(15,3.2){\makebox(0,0)[cb]{\small $2$}}
\put(16,3.25){\makebox(0,0)[cb]{\small $\tilde{2}$}}
\put(17,3.2){\makebox(0,0)[cb]{\small $2'$}}
\put(18,3.2){\makebox(0,0)[cb]{\small $\tilde{2}'$}}
\put(18.8,3.1){\makebox(0,0)[lc]{$\tilde\kappa^A_{1/2/3}$}}
\qbezier[100](14,3)(15,2.3)(16,3)
\qbezier[70](14,3)(14.7,2.5)(15.4,2.5)
\qbezier[70](18,3)(17.3,2.5)(16.6,2.5)
\put(15.4,2.5){\line(1,0){1.2}}
\qbezier[80](14,3)(14.75,2.5)(15.5,3)
\qbezier[80](18,3)(17.25,2.5)(16.5,3)
\qbezier[80](15.5,3)(16,3.33)(16.5,3)
\multiput(0,1)(1,0){4}{\circle*{0.15}}
\multiput(5,1)(1,0){4}{\circle*{0.15}}
\put(10,1){\circle*{0.15}}
\multiput(10.7,1)(0.2,0){4}{\circle*{0.03}}
\put(0,1){\makebox(0,-0.25)[ct]{\small $1\vphantom{'}$}}
\put(1,1){\makebox(0,-0.22)[ct]{\small $\tilde{1}$}}
\put(2.1,1.1){\makebox(0,0)[lb]{\small $1'$}}
\put(3.1,1.1){\makebox(0,0)[lb]{\small $\tilde{1}'$}}
\put(5,1){\makebox(0,-0.25)[ct]{\small $2\vphantom{'}$}}
\put(5.9,1){\makebox(0,-0.2)[ct]{\small $\tilde{2}$}}
\put(7.1,1){\makebox(0,-0.15)[ct]{\small $2'$}}
\put(7.9,1){\makebox(0,-0.12)[ct]{\small $\tilde{2}'$}}
\put(10,1){\makebox(0,-0.3)[ct]{\small $3_1$}}
\put(12,1.1){\makebox(0,0)[lc]{$\kappa^B_1$}}
\qbezier[80](8.5,1.8)(9.5,1.8)(10,1)
\put(3,1.8){\line(1,0){5.5}}
\qbezier[80](1.6,1)(1.6,1.8)(3,1.8)
\qbezier[60](1.6,1)(1.6,0.4)(2.6,0.4)
\qbezier[70](2.6,0.4)(3.4,0.4)(3.9,1)
\qbezier[70](3.9,1)(4.5,1.6)(5,1.6)
\qbezier[80](6,1)(6.3,0.15)(7.5,0.15)
\qbezier[80](7.5,0.15)(8.5,0.15)(8.5,1)
\qbezier[80](8.5,1)(8.5,1.6)(7.5,1.6)
\put(5,1.6){\line(1,0){2.5}}
\end{picture}
\end{center}
\caption{}
\end{figure}
\begin{figure}[ht]
\begin{center}
\setlength{\unitlength}{5mm}
\begin{picture}(22,5.6)(0,0.2)
\multiput(0,5)(1,0){4}{\circle*{0.15}}
\multiput(5,5)(1,0){4}{\circle*{0.15}}
\multiput(10,5)(2,0){6}{\circle*{0.15}}
\multiput(10.7,5)(0.2,0){4}{\circle*{0.03}}
\multiput(12.7,5)(0.2,0){4}{\circle*{0.03}}
\multiput(16.7,5)(0.2,0){4}{\circle*{0.03}}
\multiput(18.7,5)(0.2,0){4}{\circle*{0.03}}
\put(0,5){\makebox(0,-0.25)[ct]{\small $1\vphantom{'}$}}
\put(1,5){\makebox(0,-0.22)[ct]{\small $\tilde{1}$}}
\put(2.1,5.1){\makebox(0,0)[lb]{\small $1'$}}
\put(3.1,5.1){\makebox(0,0)[lb]{\small $\tilde{1}'$}}
\put(5,5){\makebox(0,-0.25)[ct]{\small $2\vphantom{'}$}}
\put(6,5){\makebox(0,-0.22)[ct]{\small $\tilde{2}$}}
\put(7.1,5){\makebox(0,-0.15)[ct]{\small $2'$}}
\put(7.9,5){\makebox(0,-0.12)[ct]{\small $\tilde{2}'$}}
\put(10,5){\makebox(0,-0.3)[ct]{\small $3_1$}}
\put(12,5){\makebox(0,-0.3)[ct]{\small $3_i$}}
\put(14,5.2){\makebox(0,0)[cb]{\small $3_{2b}$}}
\put(16,5){\makebox(0,-0.35)[ct]{\small $4_1$}}
\put(18,5){\makebox(0,-0.35)[ct]{\small $4_i$}}
\put(20,5){\makebox(0,-0.35)[ct]{\small $4_{2b}$}}
\qbezier[90](18,5)(18,5.7)(16.5,5.7)
\qbezier[90](15,5)(15,5.7)(16.5,5.7)
\qbezier[90](15,5)(15,4.5)(13.5,4.5)
\qbezier[80](12.3,5)(12.3,4.5)(13.5,4.5)
\qbezier[80](12.3,5)(12.3,6)(11,6)
\qbezier[60](12,5)(12,5.8)(11,5.8)
\put(3,5.8){\line(1,0){8}}
\put(2.8,6){\line(1,0){8.2}}
\qbezier[85](1.4,5.1)(1.4,6)(2.8,6)
\qbezier[80](1.6,5.1)(1.6,5.8)(3,5.8)
\qbezier[80](1.6,5.1)(1.6,4.6)(3,4.6)
\qbezier[100](1.4,5.1)(1.4,4.4)(3,4.4)
\qbezier[50](3,4.6)(3.5,4.6)(3.9,5)
\qbezier[60](3,4.4)(3.5,4.4)(4.1,5)
\qbezier[60](4.1,5)(4.5,5.4)(5.25,5.4)
\qbezier[60](3.9,5)(4.5,5.6)(5,5.6)
\qbezier[70](5.25,5.4)(6.3,5.4)(6.3,5)
\qbezier[80](6.3,5)(6.3,4.15)(7.5,4.15)
\qbezier[60](7.5,4.15)(8.5,4.15)(8.5,5)
\qbezier[60](8.5,5)(8.5,5.6)(7.5,5.6)
\put(5,5.6){\line(1,0){2.5}}
\put(22,5){\makebox(0,0)[lc]{$\tau'_i$}}
\multiput(0,3)(1,0){4}{\circle*{0.15}}
\multiput(5,3)(1,0){4}{\circle*{0.15}}
\multiput(10,3)(2,0){6}{\circle*{0.15}}
\multiput(10.7,3)(0.2,0){4}{\circle*{0.03}}
\multiput(12.7,3)(0.2,0){4}{\circle*{0.03}}
\multiput(16.7,3)(0.2,0){4}{\circle*{0.03}}
\multiput(18.7,3)(0.2,0){4}{\circle*{0.03}}
\put(0,3){\makebox(0,-0.3)[ct]{\small $1\vphantom{'}$}}
\put(1,3){\makebox(0,-0.27)[ct]{\small $\tilde{1}$}}
\put(2.1,3.1){\makebox(0,0)[lb]{\small $1'$}}
\put(3.1,3.1){\makebox(0,0)[lb]{\small $\tilde{1}'$}}
\put(5,3){\makebox(0,-0.3)[ct]{\small $2\vphantom{'}$}}
\put(6,3){\makebox(0,-0.27)[ct]{\small $\tilde{2}$}}
\put(7.1,3){\makebox(0,-0.15)[ct]{\small $2'$}}
\put(7.9,3){\makebox(0,-0.12)[ct]{\small $\tilde{2}'$}}
\put(10,3){\makebox(0,-0.35)[ct]{\small $3_1$}}
\put(12.1,3.1){\makebox(0,0)[lb]{\small $3_i$}}
\put(14,3.2){\makebox(0,0)[cb]{\small $3_{2b}$}}
\put(16,3){\makebox(0,-0.35)[ct]{\small $4_1$}}
\put(18,3){\makebox(0,-0.35)[ct]{\small $4_i$}}
\put(20,3){\makebox(0,-0.35)[ct]{\small $4_{2b}$}}
\qbezier[90](18,3)(18,3.7)(16.5,3.7)
\qbezier[90](15,3)(15,3.7)(16.5,3.7)
\qbezier[90](15,3)(15,2.3)(13.5,2.3)
\qbezier[100](11.7,3)(11.7,2.3)(13.5,2.3)
\qbezier[65](12,3)(12,4)(11,4)
\qbezier[50](11.7,3)(11.7,3.85)(11,3.85)
\put(3,3.85){\line(1,0){8}}
\put(2.8,4){\line(1,0){8.2}}
\qbezier[85](1.4,3.1)(1.4,4)(2.8,4)
\qbezier[80](1.6,3.1)(1.6,3.85)(3,3.85)
\qbezier[80](1.6,3.1)(1.6,2.6)(3,2.6)
\qbezier[100](1.4,3.1)(1.4,2.4)(3,2.4)
\qbezier[60](3,2.6)(3.5,2.6)(3.9,3.1)
\qbezier[70](3,2.4)(3.5,2.4)(4.1,3.1)
\qbezier[60](4.1,3.1)(4.5,3.55)(5,3.55)
\qbezier[60](3.9,3.1)(4.5,3.7)(5,3.7)
\qbezier[60](5.25,3.4)(6.3,3.4)(6.3,3)
\qbezier[70](6.3,3)(6.3,2.1)(7.5,2.1)
\qbezier[70](7.5,2.1)(8.5,2.1)(8.5,3)
\qbezier[60](8.5,3)(8.5,3.55)(7.5,3.55)
\put(5,3.7){\line(1,0){2.5}}
\put(5,3.55){\line(1,0){2.5}}
\qbezier[80](7.5,3.7)(8.7,3.7)(8.7,3)
\qbezier[80](8.7,3)(8.7,1.95)(7.5,1.95)
\qbezier[60](5.25,3.4)(4.3,3.4)(4.3,2.8)
\qbezier[70](4.3,2.8)(4.3,1.95)(5.5,1.95)
\put(5.5,1.95){\line(1,0){2}}
\put(22,3){\makebox(0,0)[lc]{$\tau''_i$}}
\multiput(0,1)(1,0){4}{\circle*{0.15}}
\multiput(5,1)(1,0){4}{\circle*{0.15}}
\put(10,1){\circle*{0.15}}
\put(14,1){\circle*{0.15}}
\put(16,1){\circle*{0.15}}
\put(20,1){\circle*{0.15}}
\multiput(10.7,1)(0.2,0){14}{\circle*{0.03}}
\multiput(16.7,1)(0.2,0){14}{\circle*{0.03}}
\put(0.1,1.1){\makebox(0,0)[lb]{\small $1\vphantom{'}$}}
\put(1.1,1.1){\makebox(0,0)[lb]{\small $\tilde{1}$}}
\put(2.35,1){\makebox(0,-0.05)[ct]{\small $1'$}}
\put(3.35,1){\makebox(0,-0.02)[ct]{\small $\tilde{1}'$}}
\put(5.3,1){\makebox(0,-0.05)[ct]{\small $2\vphantom{'}$}}
\put(6.35,1){\makebox(0,-0.02)[ct]{\small $\tilde{2}$}}
\put(7.35,1){\makebox(0,-0.05)[ct]{\small $2'$}}
\put(8.35,1){\makebox(0,-0.02)[ct]{\small $\tilde{2}'$}}
\put(10,1){\makebox(0,-0.2)[ct]{\small $3_1$}}
\put(14,1){\makebox(0,-0.2)[ct]{\small $3_{2b}$}}
\put(16,1){\makebox(0,-0.25)[ct]{\small $4_1$}}
\put(20,1){\makebox(0,-0.25)[ct]{\small $4_{2b}$}}
\put(22,1){\makebox(0,0)[lc]{$t,\tilde{t}$}}
\qbezier[80](1,1)(1.7,0.2)(2.5,0.2)
\qbezier[100](0,1)(1,0)(2.5,0)
\qbezier[80](3,1)(3.5,1.5)(4.2,1.5)
\qbezier[100](2,1)(2.8,1.7)(4.2,1.7)
\qbezier[80](13.5,1.5)(14.8,1.5)(14.8,0.9)
\qbezier[100](13.5,1.7)(15.2,1.7)(15.2,0.9)
\qbezier[50](14.8,0.9)(14.8,0.2)(14,0.2)
\qbezier[70](15.2,0.9)(15.2,0)(14,0)
\put(2.5,0.2){\line(1,0){11.5}}
\put(2.5,0){\line(1,0){11.5}}
\put(4.2,1.5){\line(1,0){9.3}}
\put(4.2,1.7){\line(1,0){9.3}}
\end{picture}
\end{center}
\caption{}
\end{figure}
Moreover, for reasons that
will be apparent later in the argument, the cusp points are also of 
limited relevance for our purposes; those which will play a role in the
argument, namely the six cusps near point A
and one of the $12b$ cusps near point B of Figure 6, give rise to braid
monodromies equal to the cubes of the half-twists represented in Figure 7.
Actually, the truly important information is contained in the vertical 
tangencies, which correspond to the half-twists
$\tau'_1,\dots,\tau'_{2b},\tau''_1,\dots,\tau''_{2b},t,\tilde{t}\in B_{4b+8}$
represented in Figure 8. 
As in \S 6.1,
the reference fiber of $\pi$ is $\{x=A\}$ for $A$ a large positive real
constant, and the chosen generating paths in the base ($x$-plane) remain 
under the real axis except near their end points; the labels 
$1,1',2,2',\tilde{1},\tilde{1}',\tilde{2},\tilde{2}'$ and
$3_1,\dots,3_{2b},4_1,\dots,4_{2b}$ correspond respectively to the two 
copies of $D_g$ and to $g(C)$.

We now turn to vertices along the right boundary of the diagram, at positions
$pj$ with $1\le j\le q-1$. The local geometric configuration is very similar
to that for the vertices along the top boundary, except for the local
description of the curve $g(C)$ which now involves $2a$ parallel copies of
$g(C_0)$ instead of $2b$. 
Another difference is that, due to the ordering of the
vertices and edges of the diagram, the slope of some of the line components
to which $g(C)$ degenerates becomes smaller than that of some of the 
components to which $D_g$
degenerates, so that the braid monodromy has to be calculated again, with
results very similar to those above. In fact, it can easily be checked that,
up to a Hurwitz equivalence, the only effect of the change of ordering on 
the local braid monodromy is the simultaneous conjugation of all contributions
by a braid that exchanges the groups of points labelled
$2,\tilde{2},2',\tilde{2}'$ and $3_1,\dots,3_{2a}$ by moving
them around each other counterclockwise.

The last vertex that remains to be investigated is the corner vertex at
position $pq$. The local configuration for $D^0=D_g\cup g(C)$ is obtained
from that represented in Figure 9 (left) by smoothing the $4ab$ mutual 
intersections between the lines labelled $2_1,\dots,2_{2a}$ and 
$3_1,\dots,3_{2b}$. Indeed, the local configuration for $D_g$ is simply
a conic (labelled $1,1'$ in Figure 9), while $g(C_0)$ consists of two
lines tangent to that conic, and $g(C)$ is obtained by ``thickening'' these
two lines into respectively $2a$ and $2b$ components ($2_1,\dots,2_{2a}$
corresponding to the vertical edge of the diagram, and $3_1,\dots,3_{2b}$
corresponding to the horizontal edge of the diagram) and smoothing their
mutual intersections. The curve $D$ is 
then obtained from $D^0$ by separating the multiplicity $2$ component $D_g$ 
into two distinct copies, while each tangent intersection of $D_g$ with
$g(C)$ gives rise to three cusps.

\begin{figure}[ht]
\begin{center}
\setlength{\unitlength}{6mm}
\begin{picture}(7,4.7)(-0.5,-0.2)
\qbezier[300](0.2,0.5)(0.9,2)(3,2.95)
\qbezier[300](0.2,0.5)(0,-0.5)(2,-0.5)
\qbezier[300](3,2.95)(4,3.45)(6,4)
\qbezier[300](2,-0.5)(4,-0.5)(6,0)
\put(6,4.2){\line(-5,-2){6.5}}
\put(6,4.47){\line(-2,-1){6.5}}
\put(5.6,4.6){\line(-5,-3){6.1}}
\put(5.25,4.6){\line(-3,-2){5.75}}
\put(0.497,1.431){\line(-3,-2){0.9}}
\put(-0.42,-0.5){\line(3,5){3}}
\put(-0.32,-0.5){\line(1,2){2.5}}
\put(-0.2,-0.5){\line(2,5){2}}
\put(-0.12,-0.5){\line(1,3){1.667}}
\put(6,0.2){\makebox(0,0)[rb]{$1$}}
\put(6,3.8){\makebox(0,0)[rt]{$1'$}}
\put(6.1,4.25){\makebox(0,0)[lc]{$2_1$}}
\put(4.77,4.5){\makebox(0,0)[rc]{$2_{2a}$}}
\put(1.47,4.5){\makebox(0,0)[rc]{$3_{2b}$}}
\put(2.65,4.5){\makebox(0,0)[lc]{$3_{1}$}}
\end{picture}
\quad\setlength{\unitlength}{5mm}
\begin{picture}(14,6)(0,0)
\multiput(0,5)(1,0){4}{\circle*{0.15}}
\multiput(4,5)(2,0){3}{\circle*{0.15}}
\multiput(9,5)(2,0){3}{\circle*{0.15}}
\multiput(4.7,5)(0.2,0){4}{\circle*{0.03}}
\multiput(6.7,5)(0.2,0){4}{\circle*{0.03}}
\multiput(9.7,5)(0.2,0){4}{\circle*{0.03}}
\multiput(11.7,5)(0.2,0){4}{\circle*{0.03}}
\put(0,5){\makebox(0,-0.2)[ct]{\small $1\vphantom{'}$}}
\put(1,5){\makebox(0,-0.17)[ct]{\small $\tilde{1}$}}
\put(2.1,5){\makebox(0,-0.12)[ct]{\small $1'$}}
\put(2.8,5){\makebox(0,-0.085)[ct]{\small $\tilde{1}'$}}
\put(4,5){\makebox(0,-0.22)[ct]{\small $2_1$}}
\put(6,5){\makebox(0,-0.22)[ct]{\small $2_i$}}
\put(8,5){\makebox(0,-0.22)[ct]{\small $2_{2a}$}}
\put(9,5){\makebox(0,-0.22)[ct]{\small $3_1$}}
\put(11,5){\makebox(0,-0.22)[ct]{\small $3_j$}}
\put(13,5){\makebox(0,-0.22)[ct]{\small $3_{2b}$}}
\qbezier[100](1.5,5)(1.5,5.8)(2.5,5.8)
\qbezier[100](1.5,5)(1.5,4.2)(2.5,4.2)
\qbezier[100](3.5,5)(3.5,4.2)(2.5,4.2)
\qbezier[100](3.5,5)(3.5,5.6)(4.75,5.6)
\qbezier[100](6,5)(6,5.6)(4.75,5.6)
\qbezier[100](11,5)(11,5.8)(9.5,5.8)
\put(2.5,5.8){\line(1,0){7}}
\put(13.7,5){\makebox(0,0)[lc]{$\tau_{ij}$}}
\multiput(0,3)(1,0){4}{\circle*{0.15}}
\multiput(4,3)(4,0){2}{\circle*{0.15}}
\multiput(9,3)(4,0){2}{\circle*{0.15}}
\multiput(4.7,3)(0.2,0){14}{\circle*{0.03}}
\multiput(9.7,3)(0.2,0){14}{\circle*{0.03}}
\put(0.1,3.1){\makebox(0,0)[lb]{\small $1\vphantom{'}$}}
\put(1.1,3.1){\makebox(0,0)[lb]{\small $\tilde{1}$}}
\put(2.35,3){\makebox(0,-0.05)[ct]{\small $1'$}}
\put(3.35,3){\makebox(0,-0.02)[ct]{\small $\tilde{1}'$}}
\put(4,3){\makebox(0,-0.2)[ct]{\small $2_1$}}
\put(8,3){\makebox(0,-0.2)[ct]{\small $2_{2a}$}}
\put(9,3){\makebox(0,-0.25)[ct]{\small $3_1$}}
\put(12.9,3.05){\makebox(0,-0.25)[ct]{\small $3_{2b}$}}
\put(13.5,2){\makebox(0,0)[lt]{$t,\tilde{t}$}}
\qbezier[100](1,3)(1.7,2.2)(2.5,2.2)
\qbezier[120](0,3)(1,2)(2.5,2)
\qbezier[80](3,3)(3.5,3.5)(4.2,3.5)
\qbezier[120](2,3)(2.8,3.7)(4.2,3.7)
\qbezier[100](12.5,3.5)(13.8,3.5)(13.8,2.9)
\qbezier[120](12.5,3.7)(14.2,3.7)(14.2,2.9)
\qbezier[100](13.8,2.9)(13.8,2.2)(13,2.2)
\qbezier[120](14.2,2.9)(14.2,2)(13,2)
\put(2.5,2.2){\line(1,0){10.5}}
\put(2.5,2){\line(1,0){10.5}}
\put(4.2,3.5){\line(1,0){8.3}}
\put(4.2,3.7){\line(1,0){8.3}}
\end{picture}
\end{center}
\caption{}
\end{figure}

The braid monodromy for the corner vertex can be deduced explicitly from
this description. We are particularly interested in 
the $8ab+2$ vertical tangencies of the local model, for which the
corresponding
half-twists $\tau_{ij}$ ($1\le i\le 2a$, $1\le j\le 2b$, each appearing twice),
$t$ and $\tilde{t}$ in $B_{2a+2b+4}$ are represented in Figure 9 (right).

\subsection{Fundamental group calculations} As in \S 6, the Zariski-Van
Kampen theorem provides an explicit presentation of $\pi_1(\C^2-D)$ in
terms of the braid monodromy. The main difference is that there are now
four generators for each interior edge of the diagram 
(Figure 5), because the regeneration process involves
two copies of the branch curve of $g$; we denote by $\gamma_i,\gamma'_i$ and
$\tilde{\gamma}_i,\tilde{\gamma}'_i$ the four generators corresponding to
the $i$-th interior edge. Moreover, each edge along the top boundary of the
diagram contributes $2b$ generators (denoted by $z_{i,1},\dots,z_{i,2b}$
for the horizontal edge in position $iq$, where $1\le i\le p$), and similarly
each edge along the right boundary contributes $2a$ generators ($y_{j,1},
\dots,y_{j,2a}$ for the vertical edge in position $pj$, where $1\le j\le q$).

We are in fact interested in the stabilized quotient $G$ of $\pi_1(\C^2-D)$
(see Definition 2.2), which can be expressed in terms of the same generators
by adding suitable commutation relations.
Let $\Gamma$ be the subgroup of $G$ generated by the
$\gamma_i,\gamma'_i$, and let $\tilde\Gamma$ be the subgroup generated
by the $\tilde\gamma_i,\tilde\gamma'_i$. By definition, the elements of
$\Gamma$ always commute with those of $\tilde\Gamma$, because
the images by the geometric monodromy representation $\theta$ of the 
geometric generators $\gamma_i,\gamma'_i$ and $\tilde\gamma_i,
\tilde\gamma'_i$ act on two disjoint sets of $n/2=2pq$ sheets of 
the branched cover $f$.

As in \S 6, we introduce twisted generators $e_i,e'_i$ and $\tilde{e}_i,
\tilde{e}'_i$ for $\Gamma$ and $\tilde\Gamma$, by choosing integers $l_i$
satisfying the same compatibility conditions at the inner vertices as in
\S 6, and setting as previously $e_i=\rho_i^{l_i}(\gamma_i)$,
$e'_i=\rho_i^{l_i}(\gamma'_i)$,
$\tilde{e}_i=\tilde\rho_i^{l_i}(\tilde\gamma_i)$ and
$\tilde{e}'_i=\tilde\rho_i^{l_i}(\tilde\gamma'_i)$, with the obvious
definition for $\rho_i$ and $\tilde\rho_i$. Even though this could be
avoided by proving a suitable invariance property, we will assume that
$l_{i}=1$ for every diagonal edge in the top-most row or in the
right-most column of the diagram
(so $e_i=\gamma'_i$, $\tilde{e}_i=\tilde\gamma'_i$), and $l_j=0$ for
every vertical edge in the top-most row and every horizontal edge in
the right-most column (so $e_j=\gamma_j$, $\tilde{e}_j=\tilde\gamma_j$).
Finally, as in \S 6.1 we let $a_i=e'_ie_i^{-1}$ and
$\tilde{a}_i=\tilde{e}'_i\tilde{e}_i^{-1}$, and we relabel these elements
as $d_{ij},v_{ij},h_{ij}$ (resp.\
$\tilde{d}_{ij},\tilde{v}_{ij},\tilde{h}_{ij}$) according to their position
in the diagram. 

\begin{lemma}
The subgroup $\mathcal{B}_\Gamma\subset\Gamma$ generated by the $e_i$
and the subgroup $\mathcal{B}_{\tilde\Gamma}\subset\tilde\Gamma$ generated
by the $\tilde{e}_i$ are naturally isomorphic to quotients of
$\tilde{B}_{n/2}$. Moreover, the subgroups $\,\Gamma$ and $\,\tilde\Gamma$
of $\,G$ are naturally isomorphic to quotients of 
$\tilde{B}_{n/2}^{\smash{(}2\smash{)}}$,
with geometric generators corresponding to pairs of half-twists.
Furthermore, $\Gamma$ is generated by the elements of $\mathcal{B}_\Gamma$
and $v_{11}$, and $\smash{\tilde\Gamma}$ is generated by the elements of
$\mathcal{B}_{\tilde\Gamma}$ and $\tilde{v}_{11}$.
\end{lemma}

\proof
We first look at relations corresponding to the interior vertices of the
diagram (Figure 5) and to the vertices along the bottom and left boundaries.
Since the local description of $D$ at these vertices simply consists of
two superimposed copies of $D_g$, and since the generators of $\Gamma$
commute with those of $\tilde\Gamma$, one easily checks that the local
configurations yield relations among the $e_i,e'_i$ that are
exactly identical to those discussed in \S 6 in the case of $\CP^1\times
\CP^1$; additionally, an identical set of relations also holds among the
$\tilde{e}_i,\tilde{e}'_i$. 

Next we consider
the local configuration at a vertex along the top boundary of
the diagram, and more precisely the cusp singularities present near the 
point labelled A on Figure 6, as pictured on Figure 7. Denoting by $i$
and $j$ respectively the labels of the diagonal and vertical edges
meeting at the given vertex, the relations
corresponding to these six cusps are
\begin{align}
&\gamma'_i\gamma_j\gamma'_i=\gamma_j\gamma'_i\gamma_j,\ %
\gamma'_i\gamma'_j\gamma'_i=\gamma'_j\gamma'_i\gamma'_j,\ %
\gamma'_i(\gamma^{-1}_j\gamma'_j\gamma_j)\gamma'_i=
(\gamma^{-1}_j\gamma'_j\gamma_j)\gamma'_i(\gamma^{-1}_j\gamma'_j\gamma_j),\\
\nonumber
&\tilde\gamma'_i\tilde\gamma_j\tilde\gamma'_i=\tilde\gamma_j\tilde\gamma'_i\tilde\gamma_j,\ %
\tilde\gamma'_i\tilde\gamma'_j\tilde\gamma'_i=\tilde\gamma'_j\tilde\gamma'_i\tilde\gamma'_j,\ %
\tilde\gamma'_i(\tilde\gamma^{-1}_j\tilde\gamma'_j\tilde\gamma_j)\tilde\gamma'_i=
(\tilde\gamma^{-1}_j\tilde\gamma'_j\tilde\gamma_j)\tilde\gamma'_i(\tilde\gamma^{-1}_j\tilde\gamma'_j\tilde\gamma_j).
\end{align}
It can easily be checked that these relations satisfy a property of
invariance under twisting similar to that of 3-points. In fact, replacing
the various generators by their images under arbitrary powers of the
twisting actions $\rho_i,\tilde\rho_i,\rho_j,\tilde\rho_j$ amounts to a
conjugation of the relations (7.1) by braids belonging to the local monodromy
(either the entire local monodromy, or two of the six cusps near A, or 
combinations thereof), and thus always yields valid relations.

Therefore, the twisted generators $e_i,e'_i,e_j,e'_j$ of $\Gamma$ satisfy 
the relations (6.2), and similarly for $\tilde{e}_i,\tilde{e}'_i,\tilde{e}_j,
\tilde{e}'_j$ in $\tilde\Gamma$. One easily checks
that a similar conclusion holds for pairs of inner edges meeting at 
a vertex along the right boundary of the diagram (recall that
the local braid monodromy only differs by a simple conjugation).
Finally, because we are looking at the stabilized fundamental group, the
commutation relations discussed in \S 6 automatically hold in $\Gamma$
and $\tilde\Gamma$. 

So, except for the equality relations arising from vertical tangencies at 
the vertices along the top and right boundaries of the diagram,
all the relations described in \S 6.1 for the case 
of $\CP^1\times\CP^1$ simultaneously hold in
$\Gamma$ and in $\tilde\Gamma$.
Therefore, the structure of $\Gamma$ and $\tilde\Gamma$ can be studied
by the same argument as in the case of $\CP^1\times\CP^1$ 
(\cite{MSegre}, see also \S 6), which yields the desired result.
\endproof

\begin{lemma}
The equality $z_{r,i}=z_{r,1}$ holds for every $1\le r\le p$, $1\le i\le
2b$; similarly, $y_{r,i}=y_{r,1}$ for every $1\le r\le q$, $1\le i\le 2a$.
Moreover, the $y_{r,i}$ and the $z_{r,i}$ are all
conjugates of $y_{q,1}$ under the action of elements of $\mathcal{B}_\Gamma$
and $\mathcal{B}_{\tilde\Gamma}$.
\end{lemma}

\proof
First consider the corner vertex at position $pq$, and more
precisely the half-twists $\tau_{ij}$ arising from the vertical tangencies
of the local model near this vertex (Figure 9). Denoting by $\mu$ the label 
of the diagonal edge in position $pq$, the half-twist $\tau_{1i}$ yields
the relation $(y_{q,1}^{-1}\dots y_{q,2a}^{-1} z_{p,1}^{-1}\dots
z_{p,i-1}^{-1})z_{p,i}(z_{p,i-1}\dots z_{p,1}y_{q,2a}\dots y_{q,1})=
\tilde\gamma'_\mu\gamma'_\mu y_{q,1}\gamma'_\mu{}^{-1}
\tilde\gamma'_\mu{}^{-1}$. It follows that the quantity 
$(z_{p,1}^{-1}\dots z_{p,i-1}^{-1})
z_{p,i}(z_{p,i-1}\dots z_{p,1})$ is independent of $i$, which by an easy
induction on $i$ implies that $z_{p,i}=z_{p,1}$ for all $i$.
Observing that $y_{q,1},\dots,y_{q,2a}$ and
$z_{p,1},\dots,z_{p,2b}$ are mapped by $\theta$ to disjoint transpositions
and hence commute in $G$, we in fact have
$z_{p,i}=\tilde\gamma'_\mu\gamma'_\mu y_{q,1}\gamma'_\mu{}^{-1}
\tilde\gamma'_\mu{}^{-1}$ for all $i$. Since by assumption the twisting
parameter $l_\mu$ is equal to $1$, the generators $\gamma'_\mu=e_\mu$ and
$\tilde\gamma'_\mu=\tilde{e}_\mu$ belong to $\mathcal{B}_\Gamma$ and 
$\mathcal{B}_{\tilde\Gamma}$ respectively. This proves the claims made about
the $z_{p,i}$. 

Similarly comparing the relations corresponding to the
half-twists $\tau_{i1}$, it can be seen immediately that the quantity
$(y_{q,1}^{-1}\dots y_{q,i-1}^{-1})y_{q,i}(y_{q,i-1}\dots y_{q,1})$ is
independent of $i$, which implies that $y_{q,i}=y_{q,1}$ for all $i$.

We now proceed by induction~: assume that $z_{r+1,i}=z_{r+1,1}$ for all $i$,
and that $z_{r+1,1}$ is a conjugate of $y_{q,1}$ under the action of
$\mathcal{B}_\Gamma$ and $\mathcal{B}_{\tilde\Gamma}$. Let $\mu$ and $\nu$ 
be the labels of the diagonal and vertical edges meeting at the
vertex in position $rq$, and let $\psi_r=\tilde\gamma'_\nu\gamma'_\nu
\tilde\gamma^{\vphantom{1}}_\nu\gamma^{\vphantom{1}}_\nu\tilde\gamma'_\mu\gamma'_\mu\gamma_\nu^{-1}
\tilde\gamma_\nu^{-1}\gamma'_\nu {}^{-1}\tilde\gamma'_\nu {}^{-1}$.
Define $\zeta_r=\psi_r\gamma_\nu\psi_r^{-1}$,
$\zeta'_r=\psi_r\gamma'_\nu\psi_r^{-1}$,
$\tilde\zeta_r=\psi_r\tilde\gamma_\nu\psi_r^{-1}$, and
$\tilde\zeta'_r=\psi_r\tilde\gamma'_\nu\psi_r^{-1}$. Recalling that the
elements of $\Gamma$ commute with those of $\tilde\Gamma$, the relations
(7.1) imply that $\zeta'_r=\gamma'_\nu\gamma_\nu\gamma'_\mu(\gamma^{-1}_\nu
\gamma'_\nu\gamma_\nu)\gamma'_\mu {}^{-1}\gamma_\nu^{-1}\gamma'_\nu{}^{-1}=
\gamma'_\nu\gamma_\nu (\gamma^{-1}_\nu
\gamma'_\nu\gamma_\nu)^{-1}\gamma'_\mu(\gamma^{-1}_\nu
\gamma'_\nu\gamma_\nu)\gamma_\nu^{-1}\gamma'_\nu{}^{-1}=
\gamma_\nu\gamma'_\mu\gamma^{-1}_\nu=\gamma'_\mu{}^{-1}\gamma_\nu\gamma'_\mu$.
Similar calculations for the other elements yield that
\begin{equation}
\zeta_r=\gamma'_\mu{}^{-1}(\gamma_\nu^{-1}\gamma'_\nu\gamma_\nu)\gamma'_\mu,
\ \tilde\zeta_r=\tilde\gamma'_\mu{}^{-1}(\tilde\gamma_\nu^{-1}
\tilde\gamma'_\nu\tilde\gamma_\nu)\tilde\gamma'_\mu,\ %
\zeta'_r=\gamma'_\mu{}^{-1}\gamma_\nu\gamma'_\mu,
\ \tilde\zeta'_r=\tilde\gamma'_\mu{}^{-1}\tilde\gamma_\nu\tilde\gamma'_\mu.
\end{equation}
Due to the choice of twisting parameters $l_\mu=1$
and $l_\nu=0$, $\zeta'_r\in\mathcal{B}_\Gamma$ and $\tilde\zeta'_r\in
\mathcal{B}_{\tilde\Gamma}$.

Since the $z_{r,i}$ commute with the $z_{r+1,i}$ in $G$ 
(they are mapped to disjoint transpositions by $\theta$), and since
by assumption $z_{r+1,i}=z_{r+1,1}$ for all $i$, we have
$$(z_{r,1}^{-1}\dots z_{r,i}^{-1}z_{r+1,1}^{-1}\dots z_{r+1,i-1}^{-1})z_{r+1,i}
(z_{r+1,i-1}\dots z_{r+1,1}z_{r,i}\dots z_{r,1})=z_{r+1,1}$$ for all $i$. 
Therefore, the relation
arising from the vertical tangency $\tau'_i$ (Figure 8) at the vertex $rq$
can be written in the form $$z_{r+1,1}=\tilde\zeta'_r\zeta'_r
(z_{r,1}^{-1}\dots z_{r,i-1}^{-1})z_{r,i}(z_{r,i-1}\dots z_{r,1})
\zeta'_r{}^{-1}\tilde\zeta'_r{}^{-1}.$$
In particular, the value of $(z_{r,1}^{-1}\dots z_{r,i-1}^{-1})z_{r,i}
(z_{r,i-1}\dots z_{r,1})$ does not depend on $i$, which implies that
$z_{r,i}=z_{r,1}$ for all $i$. Moreover, we have
$z_{r,i}=\zeta'_r{}^{-1}\tilde\zeta'_r{}^{-1} z_{r+1,1}\tilde\zeta'_r
\zeta'_r.$
So, by induction on decreasing values of $r$, we obtain the desired
results about $z_{r,i}$. The case of $y_{r,i}$ is handled using exactly the
same argument, going inductively through the vertices along the
right boundary of the diagram. Indeed, observe that the local braid monodromy
at one of these vertices simply differs from that at a vertex along the top
boundary by a conjugation which exchanges the positions of two groups of
geometric generators~; however, because the corresponding transpositions in
$S_n$ are disjoint, these generators commute with each other in $G$, so that
the relations induced by the local braid monodromy can be expressed in
exactly the same form.
\endproof

\begin{lemma}
The element $\tilde{v}_{11}$ belongs to the subgroup of $G$ generated by
$\Gamma$, $\mathcal{B}_{\tilde\Gamma}$, and $y_{q,1}$.
\end{lemma}

\proof
Consider the local relations for the vertex at position $1q$, and more
precisely the equality relation corresponding to the half-twist labelled
$\tau''_1$ in Figure 8~: with the same notations as in the proof of Lemma
7.2, we have $z_{2,1}=\zeta_1^{-1}\tilde\zeta_1^{-1}z_{1,1}\tilde\zeta_1
\zeta_1$. Moreover, the cusp point with monodromy $\kappa_1^{B}$ pictured
on Figure 7 yields the relation $\tilde\zeta_1 z_{1,1} \tilde\zeta_1=
z_{1,1} \tilde\zeta_1 z_{1,1}$. It follows that $z_{2,1}=\zeta_1^{-1}
z_{1,1}\tilde\zeta_1 z_{1,1}^{-1} \zeta_1$. Therefore, using formula (7.2)
for $\tilde\zeta_1$, we obtain $\tilde\gamma'_\nu=\tilde\gamma_\nu
\tilde\gamma'_\mu z_{1,1}^{-1}\zeta_1 z_{2,1}
\zeta_1^{-1}z_{1,1}\tilde\gamma'_\mu{}^{-1}\tilde\gamma_\nu^{-1}$, 
where $\mu$ and $\nu$ are the labels of the two interior
edges meeting at the considered vertex. 

Observe that, since $l_\nu=0$ and $l_\mu=1$, the generators
$\tilde\gamma_\nu=\tilde{e}_\nu$ and $\tilde\gamma'_\mu=
\tilde{e}_\mu$ belong to $\mathcal{B}_{\tilde\Gamma}$. Moreover, it is obvious
from (7.2) that $\zeta_1\in \Gamma$. Using the result of Lemma 7.2 to
express $z_{1,1}$ and $z_{2,1}$ in terms of $y_{q,1}$, it follows that
$\tilde\gamma'_\nu=\tilde{e}'_\nu$ belongs to the subgroup of $G$ generated
by $\Gamma$, $\mathcal{B}_{\tilde\Gamma}$, and $y_{q,1}$. Therefore,
$\tilde{v}_{1,q}=\tilde{e}'_\nu\tilde{e}_\nu^{-1}$ also belongs to this
subgroup. Finally, the local relations analogous to (6.5) 
for the $\tilde{e}_i$ and $\tilde{a}_i$ at the vertex in position $1r$ imply
that $\tilde{v}_{1,r}$ and $\tilde{v}_{1,r+1}$ are conjugates of each
other under
the action of elements of $\mathcal{B}_{\tilde\Gamma}$. Therefore, by
induction $\tilde{v}_{1,1}$ can be expressed in terms of $\tilde{v}_{1,q}$
and elements of $\mathcal{B}_{\tilde\Gamma}$, which completes the proof.
\endproof

\begin{lemma}
The subgroup $\mathcal{B}$ of $G$ generated by $\mathcal{B}_\Gamma$,
$\mathcal{B}_{\tilde\Gamma}$ and $y_{q,1}$ is naturally a quotient
of $\tilde{B}_n$, with geometric generators corresponding to half-twists.
\end{lemma}

\proof
We construct a surjective map $\alpha:\tilde{B}_n\to\mathcal{B}$ as follows
(recall that $n=4pq$).
First observe that the subgroup of $\tilde{B}_n$ generated by the
half-twists $x_1,\dots,x_{2pq-1}$ is naturally isomorphic to
$\tilde{B}_{n/2}$, which by Lemma 7.1 admits a surjective homomorphism to
$\mathcal{B}_\Gamma$ mapping half-twists to geometric generators.
We use this homomorphism to define $\alpha(x_i)$ for $1\le i\le 2pq-1$.
Any two half-twists in $\tilde{B}_{n/2}$ are conjugate to each 
other; therefore, after a suitable conjugation we can assume that
$\alpha(x_{2pq-1})=e_{\mu}$, where $\mu$ is the label of the diagonal edge 
at position $pq$ in the diagram, and that the other $\alpha(x_i)$
($i\le 2pq-2$) are geometric generators mapped by $\theta$ to transpositions
disjoint from $\theta(y_{q,1})$. Because of the stabilization process,
this last requirement implies that $\alpha(x_i)$ commutes with $y_{q,1}$
for $i\le 2pq-2$.

Similarly, the subgroup of $\tilde{B}_n$
generated by $x_{2pq+1},\dots,x_{n-1}$ is naturally isomorphic to
$\tilde{B}_{n/2}$ and admits a surjective homomorphism to
$\mathcal{B}_{\tilde\Gamma}$, which we use to define $\alpha(x_i)$ for
$2pq+1\le i\le n-1$. Once again, without loss of generality we can assume
that $\alpha(x_{2pq+1})=\tilde{e}_\mu$ and that the other $\alpha(x_i)$
commute with $y_{q,1}$. Finally, we define 
$\alpha(x_{2pq})=y_{q,1}$.

All that remains to be checked is that $\alpha$ can be made into a group
homomorphism (obviously surjective by construction), i.e.\ that the
relations defining $\tilde{B}_n$ are also satisfied by the chosen images
$\alpha(x_i)$ in $\mathcal{B}$. Since $\alpha$ is built out of two group
homomorphisms and since the elements of $\mathcal{B}_\Gamma$ commute with
those of $\mathcal{B}_{\tilde\Gamma}$, the only relations to be checked
are those involving $x_{2pq}$.

Consider the corner vertex at position $pq$ in the diagram:
the cusp singularities arising from the regeneration of the rightmost 
tangent intersection of $D_g$ with $g(C)$ in Figure 9 imply the relations 
$\gamma'_\mu y_{q,1}\gamma'_\mu=y_{q,1}\gamma'_\mu y_{q,1}$ and 
$\tilde\gamma'_\mu y_{q,1}\tilde\gamma'_\mu =y_{q,1}\tilde\gamma'_\mu
y_{q,1}$. Since $l_\mu=1$, we have $\gamma'_\mu=e_\mu$ and
$\tilde\gamma'_\mu=\tilde{e}_\mu$, so that these relations can be
rewritten as $\alpha(x_{2pq-1})\alpha(x_{2pq})\alpha(x_{2pq-1})=
\alpha(x_{2pq})\alpha(x_{2pq-1})\alpha(x_{2pq})$ and
$\alpha(x_{2pq+1})\alpha(x_{2pq})\alpha(x_{2pq+1})=
\alpha(x_{2pq})\alpha(x_{2pq+1})\alpha(x_{2pq})$. Finally, for all
$i$ such that $|i-2pq|\ge 2$, the relation $[\alpha(x_{2pq}),\alpha(x_i)]=1$
holds by construction.
Therefore, $\alpha$ 
defines a surjective group homomorphism from $\tilde{B}_n$ to $\mathcal{B}$,
mapping half-twists to geometric generators.
\endproof

\begin{proposition}
The morphism $\alpha$ extends to a surjective group homomorphism
from $\tilde{B}_n^{(2)}\simeq \tilde{B}_n \ltimes \tilde{P}_{n,0}$ to
$G$ mapping pairs of half-twists to geometric generators. In particular,
the group $G$ has property $(*)$.
\end{proposition}

\proof
Lemma 7.2 implies that $G$ is generated by $\Gamma$, $\tilde\Gamma$, and
$y_{q,1}$. Therefore, by Lemma 7.1, $G$ is generated by $\mathcal{B}$,
$v_{11}$ and $\tilde{v}_{11}$, while Lemma 7.3 implies that $\tilde{v}_{11}$
can be eliminated from the list of generators. Since Lemma 7.4 identifies
$\mathcal{B}$ with a quotient of $\tilde{B}_n$, the main remaining task
is to check that the subgroup $\mathcal{P}$ generated by the $g^{-1}v_{11}g$,
$g\in\mathcal{B}$, is naturally isomorphic to a quotient of
$\tilde{P}_{n,0}$. This can be done by proving that $\mathcal{P}$ is
a {\em primitive $\tilde{B}_n$-group} (Definition 5 of \cite{MSegre}), as
it follows from the discussion in \S 1 of \cite{MSegre} that every such
group is a quotient of $\tilde{P}_{n,0}$ (compare Propositions 1, 2, 3 of
\cite{MSegre} with the presentation of $\tilde{P}_{n,0}$ given in Lemma
3.1). 

As stated in Lemma 7.1, the arguments of \cite{MSegre} show that
the subgroup generated by the $g^{-1}v_{11}g$, $g\in \mathcal{B}_\Gamma$,
is a primitive $\tilde{B}_{n/2}$-group (and hence a quotient of
$\tilde{P}_{n/2,0}$). The desired result about
$\mathcal{P}$ then follows simply by observing that $v_{11}$ commutes
with $y_{q,1}$ and with the generators of $\mathcal{B}_{\tilde\Gamma}$
and using a criterion due to Moishezon (Proposition 6 of \cite{MSegre});
indeed, an obvious corollary of this criterion is that, upon enlarging
the conjugation action from $\tilde{B}_{n/2}$ to $\tilde{B}_n$, it is
sufficient to check that the additional half-twist generators act trivially
on the given prime element ($v_{11}$).

Since $G$ is obviously generated by its subgroups $\mathcal{B}$ and
$\mathcal{P}$, and since $\mathcal{P}$ is normal, it is naturally a
quotient of $\tilde{B}_n\ltimes \tilde{P}_{n,0}\simeq \tilde{B}_n^{(2)}$.
Moreover, the geometric generators of $G$ are all mutually conjugate
(because the curve $D$ is irreducible), and by construction the $e_i$ (and
$\tilde{e}_i$) correspond to pairs of half-twists in $\tilde{B}_n^{(2)}$, 
so the same is true of all geometric generators. Finally,
by going carefully over the construction, it is not hard to check that
the end points of the half-twists $(x,y)$ corresponding to a given geometric
generator $\gamma$ are always the natural ones, in the sense that
$\sigma(x)=\sigma(y)=\theta(\gamma)$. Therefore, $G$ has property $(*)$.
\endproof

At this point, the only remaining task in the proof of Theorem 4.6 is to 
characterize the kernel of the surjective morphism
$\alpha:\tilde{B}_n^{(2)}\to G$ given by Proposition 7.5.
As a consequence of Lemmas 3.3 and 3.4, the commutation relations induced
either by nodes in the branch curve $D$ or by the stabilization process,
as well as the relations induced by the cusp points of $D$, automatically
hold, so that $\mathrm{Ker}\,\alpha$ is generated by equality relations
between pairs of half-twists induced by the vertical tangencies of $D$.
Moreover, as in \S 6.1.2 the classification of half-twists in $\tilde{B}_n$ 
(Lemma 3.2) allows us to associate to every $a_i$ (resp.\ $\tilde{a}_i$) a
pair of integers $\bar{a}_i$ (resp.\ $\tilde{\bar{a}}_i$), well-defined 
modulo the subgroup
$\Lambda=\{(\kappa,\lambda),\ (u_1^\kappa\eta^{\kappa(\kappa-1)/2},
u_1^\lambda\eta^{\lambda(\lambda-1)/2})\in\mathrm{Ker}\,\alpha\}\subset
\Z^2$. Recall however from \S 6.1.2 that this construction requires us to 
choose an ordering of the $n=4pq$ sheets of the branched cover; in our
case, these split into two sets of $2pq$ sheets, the first one on which
the $\theta(e_i),\theta(e'_i)$ act by permutations, and the second one on
which the $\theta(\tilde{e}_i),\theta(\tilde{e}'_i)$ act by permutations.
The ordering we will consider is obtained by enumerating first the first
set of $2pq$ sheets, and then the second one. In each set, the sheets are
naturally in correspondence with the $2pq$ triangles of the diagram in
Figure 5: the ordering we choose for each of the two sets of $2pq$ sheets
is obtained as in the case of $\CP^1\times\CP^1$ \cite{MSegre} by enumerating
the $2pq$ triangles of the diagram from left to right and from bottom to
top.

We have seen above that the relations coming from the vertical tangencies at
the inner vertices of the diagram and at those along the lower and left 
boundaries are exactly the same as in the case of $\CP^1\times\CP^1$,
except they simultaneously apply to the generators of $\Gamma$ 
and to those of $\tilde\Gamma$. Therefore, as in \S 6.1.2, these relations
do not contribute to $\mathrm{Ker}\,\alpha$ by themselves, but they
translate into equalities between the $\bar{a}_i$ (and similarly 
between the $\tilde{\bar{a}}_i$), which yield the following formulas 
(with the obvious notations)~: $\bar{d}_{i,j}=\,\tilde{\!\!\bar{d}}_{i,j}=
(j-i,0)$, $\bar{v}_{i,j}=\tilde{\bar{v}}_{i,j}=(1-i,1)$,
$\bar{h}_{i,j}=\tilde{\bar{h}}_{i,j}=(1-j,1)$ (compare with (6.9)).

Next, we consider the corner vertex at position $pq$, for which the
braid monodromy contribution of the vertical tangencies is represented
in Figure 9. Recall that some of the half-twists $\tau_{ij}$ were used
in the proof of Lemma 7.2 to eliminate $y_{q,2},\dots,y_{q,2a}$ and
$z_{p,1},\dots,z_{p,2b}$ from the list of generators by expressing them
in terms of $y_{q,1}$; however, since these relations imply that 
$y_{q,i}=y_{q,1}$ and $z_{p,i}=z_{p,1}$ (cf.\ Lemma 7.2), all the other
relations coming from the $\tau_{ij}$ become redundant. Therefore these
equality relations do not make any contributions to the kernel of $\alpha$.
We are left with the two half-twists $t,\tilde{t}$ of Figure 9. Denote by
$\mu$ the label of the diagonal edge passing through the corner vertex.
Because $G$ has property $(*)$, and using the results of \S 3, we can 
find an element $g\in\tilde{B}_n^{\smash{(2)}}$ such that 
$z_{p,1}=\alpha(g^{-1}(x_1,x_1)g)$,
$e_\mu=\gamma'_\mu=\alpha(g^{-1}(x_2,x_2)g)$, and $y_{q,1}=\alpha(g^{-1}
(x_3,x_3)g)$.
Recalling that $\bar{d}_{p,q}=(q-p,0)$ and observing that the conjugation
by $g$ preserves the ordering of the end points for $e_\mu$, by definition
of $\bar{d}_{p,q}$ we have $e'_\mu=\alpha(g^{-1}(x_2 u_2^{p-q}
\eta^{(p-q)(p-q-1)/2},x_2)g)$, and therefore $\gamma_\mu=e_\mu^{-1} e'_\mu
e_\mu=\alpha(g^{-1}(x_2 u_2^{q-p} \eta^{(q-p)(q-p-1)/2},x_2)g)$.
The half-twist $t$ yields the relation $\gamma_\mu=z_{p,1}^{2b} y_{q,1}^{2a}
\gamma'_\mu y_{q,1}^{-2a} z_{p,1}^{-2b}$; an easy computation shows that
the right-hand side of this relation is equal to 
$\alpha(g^{-1}(x_2 u_2^{a-b}\eta^{(a-b)(a-b-1)/2},
x_2 u_2^{a-b}\eta^{(a-b)(a-b-1)/2})g)$. Comparing the two formulas for
$\gamma_\mu$, we conclude that the relation introduced by the half-twist
$t$ is equivalent to the property that $(a-b+p-q,a-b)\in\Lambda$.
A similar calculation shows that the relation introduced by $\tilde{t}$
can also be rewritten in the form $(a-b+p-q,a-b)\in\Lambda$.

We now consider the vertex at position $rq$ ($1\le r\le p-1$), and
investigate in the same manner the equality relations coming from the
vertical tangencies $\tau'_i,\tau''_i,t,\tilde{t}$ represented in Figure 8.
Recall that the relations induced by $\tau'_i$ were used in the proof of
Lemma 7.2 to show that
$z_{r,i}=\zeta'_r{}^{-1}\tilde\zeta'_r{}^{-1}z_{r+1,1}\tilde\zeta'_r
\zeta'_r$ and consequently eliminate the $z_{r,i}$ from the list of 
generators; these relations are therefore already accounted for. Next, 
we turn to the relation induced by $\tau''_i$, which taking into account 
that $z_{r,i}=z_{r,1}$ and $z_{r+1,i}=z_{r+1,1}$ can be written in the
form $z_{r+1,1}=\zeta_r^{-1}\tilde\zeta_r^{-1}z_{r,1}\tilde\zeta_r\zeta_r$.
Using the expression of $z_{r,1}$ in terms of $z_{r+1,1}$, this identity
can also be expressed by the commutation relation $[z_{r+1,1},
\tilde\zeta'_r \zeta'_r\tilde\zeta_r\zeta_r]=1$. By (7.2),
we have $\tilde\zeta'_r \zeta'_r\tilde\zeta_r\zeta_r=e_\mu^{-1}
\tilde{e}_\mu^{-1}\tilde{e}'_\nu\tilde{e}_\nu e'_\nu e_\nu\tilde{e}_\mu
e_\mu$, where $\mu$ and $\nu$ are
the labels of the two interior edges meeting at position $rq$. 
Since $z_{r+1,1}$ commutes with $e_\mu$ and $\tilde{e}_\mu$, 
the relation can then be rewritten as $[z_{r+1,1},\tilde{e}'_\nu\tilde{e}_\nu 
e'_\nu e_\nu]=1$. Taking into account the ordering of the sheets of the
branched cover, an easy calculation in $\tilde{B}_n^{(2)}$ shows that
this relation automatically holds as a consequence of the equality
$\bar{v}_{r,q}=\tilde{\bar{v}}_{r,q}$.

The relation induced by the half-twist $t$ (Figure 8) can be expressed as
$\gamma_\mu=z_{r,1}^{2b}\gamma'_\nu \gamma_\nu \gamma'_\mu \gamma_\nu^{-1}
\gamma'_\nu{}^{-1} z_{r,1}^{-2b}$. Using property $(*)$ and recalling that
$\bar{d}_{r,q}=(q-r,0)$ and $\bar{v}_{r,q}=(1-r,1)$, we can find $g\in
\tilde{B}_n^{(2)}$, preserving the ordering of the end points for $e_\mu$
and $e_\nu$, such that $z_{r,1}=\alpha(g^{-1}(x_1,x_1)g)$, 
$\gamma'_\mu=e_\mu=\alpha(g^{-1}(x_2,x_2)g)$, $\gamma_\mu=e_\mu^{-1}
e'_\mu e_\mu=\alpha(g^{-1}(x_2 u_2^{q-r} \eta^{(q-r)(q-r-1)/2},x_2)g)$,
$\gamma_\nu=e_\nu=\alpha(g^{-1}(x_3,x_3)g)$, and $\gamma'_\nu=e'_\nu=
\alpha(g^{-1}(x_3 u_3^{r-1} \eta^{(r-1)(r-2)/2},x_3 u_3^{-1} \eta)g)$.
So $z_{r,1}^{2b}\gamma'_\nu \gamma_\nu \gamma'_\mu \gamma_\nu^{-1}
\gamma'_\nu{}^{-1} z_{r,1}^{-2b}$ is equal to
$\alpha(g^{-1}(x_2 u_2^{2-r-b}\eta^{(2-r-b)(1-r-b)/2},x_2 u_2^{2-b}
\eta^{(2-b)(1-b)/2})g)$. Comparing this with the expression for
$\gamma_\mu$, it becomes apparent that the relation induced by $t$ 
is in fact equivalent to the condition $(q+b-2,b-2)\in\Lambda$.
A similar calculation for the half-twist $\tilde{t}$ shows that the
relation it induces can also be expressed in the form
$(q+b-2,b-2)\in\Lambda$.

Finally, the case of the vertices along the right boundary of the
diagram can be studied by exactly the same argument; the relations 
corresponding to the vertical tangencies of the local model can
be expressed by the single requirement that $(p+a-2,a-2)\in\Lambda$.

Therefore, $\Lambda\subset\Z^2$ is the subgroup generated by
$(p+a-2,a-2)$ and $(q+b-2,b-2)$, and $\mathrm{Ker}\,\alpha$ is the
normal subgroup of $\tilde{B}_n^{(2)}$ generated by the two elements
$g_1=(u_1^{p+a-2}\eta^{\lambda(p+a-2)},u_1^{a-2}\eta^{\lambda(a-2)})$ and
$g_2=(u_1^{q+b-2}\eta^{\lambda(q+b-2)},u_1^{b-2}\eta^{\lambda(b-2)})$, where
$\lambda(i)=i(i-1)/2$.
Observe that $G^0_{p,q}=(\tilde{P}_{n,0}\times\tilde{P}_{n,0})/
\mathrm{Ker}\,\alpha$, and recall from Lemma 3.1 that
$[\tilde{P}_{n,0},\tilde{P}_{n,0}]=\{1,\eta\}\simeq \Z_2$ and
$\mathrm{Ab}\,\tilde{P}_{n,0}\simeq \Z^{n-1}$.

We first consider the commutator subgroup
$[G^0_{p,q},G^0_{p,q}]\simeq C/(C\cap \mathrm{Ker}\,\alpha)$, where
$C=\{1,\eta\}\times \{1,\eta\}$. First of all,
if $a+p$ is odd, then considering the commutator of
$g_1$ with
$(u_2,1)$ we obtain that $(\eta,1)\in \mathrm{Ker}\,\alpha$, and
similarly if $b+q$ is odd; otherwise, one easily checks that
$(\eta,1)\not\in\mathrm{Ker}\,\alpha$.
Moreover, if $a$ is odd, then considering the commutator of $g_1$ with
$(1,u_2)$ we obtain that $(1,\eta)\in \mathrm{Ker}\,\alpha$, and similarly
if $b$ is odd; when $a$ and $b$ are both even,
$(1,\eta)\not\in\mathrm{Ker}\,\alpha$. Also, it is easy to check that
$\mathrm{Ker}\,\alpha$ only contains $(\eta,\eta)$ if it also contains
$(\eta,1)$ and $(1,\eta)$.
The claim made in the statement of Theorem 4.6 about the structure
of $[G^0_{p,q},G^0_{p,q}]$ follows.

Finally, we have $\mathrm{Ab}\,G^0_{p,q}\simeq (\tilde{P}_{n,0}\times
\tilde{P}_{n,0})/\langle C,\mathrm{Ker}\,\alpha\rangle\simeq
(\Z^2/\Lambda)^{n-1}$. Observing that $\Z^2/\Lambda=\Z^2/\langle
(p+a-2,a-2),(q+b-2,b-2)\rangle\simeq \Z^2/\langle (p,a-2),(q,b-2)\rangle$,
this completes the proof of Theorem 4.6.


\subsection*{Acknowledgements} The authors wish to thank 
M. Gromov for prompting their interest in this problem and 
for useful advice and suggestions. The first and third authors wish to 
thank respectively Ecole Polytechnique and Imperial College for their
hospitality.

\end{document}